\documentclass[10pt]{article}
\usepackage{amsmath, amssymb}
\usepackage{latexsym}
\usepackage{mathtools}
\usepackage{ulem}
\usepackage{multirow}
\usepackage{graphicx}
\usepackage{color}
\usepackage{url}
\usepackage{framed}
\usepackage{float}
\usepackage{enumitem}
\usepackage{lipsum}
\usepackage{enumitem}
\usepackage{tabularx}
\newcommand{\subscript}[2]{$#1 _ #2$}
\providecommand{\keywords}[1]{\textbf{Keywords:} #1}

\numberwithin{equation}{section}
\DeclareMathAlphabet{\itbf}{OML}{cmm}{b}{it}

\newcommand{\RR}{\mathbb{R}}

\newcommand{\NN}{\mathbb{N}}

\newcommand{\ds}{\displaystyle}
\newcommand{\no}{\nonumber}
\newcommand{\ri}{\rightarrow}
\newcommand{\sm}{\setminus}

\newcommand{\q}{\quad}

\newcommand{\bm}{{\itbf m}}

\newcommand{\bq}{{\itbf q}}

\newcommand{\bx}{{\itbf x}}

\newcommand{\bn}{{\itbf n}}

\newcommand{\bi}{\begin{itemize}}
\newcommand{\ei}{\end{itemize}}

\newcommand{\ca}{{\cal A}}
\newcommand{\cp}{{\cal P}}

\newcommand{\bH}{{\itbf H}}

\newcommand{\be}{\begin{eqnarray}}
\newcommand{\ee}{\end{eqnarray}}

\newcommand{\ben}{\begin{eqnarray*}}
\newcommand{\een}{\end{eqnarray*}}
\def\ds{\displaystyle}

\newcommand\ov{\overline}

\newtheorem{lem}{Lemma}[section]
\newtheorem{prop}{Proposition}[section]

\newtheorem{thm}{Theorem}[section]

\newcommand{\bea}{\begin{eqnarray*}}
\newcommand{\eea}{\end{eqnarray*}}
\newcommand{\bean}{\begin{eqnarray}}
\newcommand{\eean}{\end{eqnarray}}
\newcommand{\p}{\partial}
\newcommand{\f}{\frac}
\newcommand{\s}{\sqrt}
\newcommand{\di}{\mbox{div }}
\newcommand{\aaa}{\mbox{$[$}}
\newcommand{\bbb}{\mbox{$]$}}

\begin{document}

% in this version the tangential vector field  are more 
% explicitly defined
\title{Stability properties for a class of inverse problems}
% and application to neural networks

\author{  Darko
Volkov \thanks{Department of Mathematical Sciences,
Worcester Polytechnic Institute, Worcester, MA 01609.
}  }
%joancs@cttc.upc.edu
%Heat and Mass Technological Center (CTTC). Technical University
%of Catalonia (UPC), Colom 11, 08222 Terrassa (Barcelona), Spain

\maketitle

%Hi Darko,
%Thank you for coming today and for your feedback!
%Here is the link to the recording: ​
% https://wpi0-my.sharepoint.com/personal/mmalone_wpi_edu/_layouts/15/onedrive.aspx?id=%2Fpersonal%2Fmmalone%5Fwpi%5Fedu%2FDocuments%2F2021%2D12%2D10%2011%2E10%2E09%20Elisa%20Negrini%2D%20PhD%20Proposal

\begin{abstract}
We establish Lipschitz  stability properties for a class of inverse problems.
In that class, the associated direct problem is formulated 
by an integral operator $\ca_m$ depending non-linearly on a parameter $m$
and operating on a function $u$. 
In the inversion step both $u$ and $m$ are unknown but we are only 
interested in recovering $m$. We discuss examples of such inverse problems
for the  elasticity equation with applications to seismology and for the 
inverse scattering problem in electromagnetic theory. 
Assuming a few injectivity and regularity properties  for $\ca_m$,
we prove that the  
 inverse problem 
with a 
 a finite number of data points
is solvable and that the solution is 
Lipschitz stable in the data.
We show a reconstruction example illustrating the use of neural networks. 
\end{abstract}
\textbf{MSC 2010 Mathematics Subject Classification:} 
% 35R30: Inverse problems (undetermined coefficients, etc.) for PDE
% 45-xx: Integral equations45Q05: Inverse problems
35R30, 47N20, 47N40.
% 45Q05 not here
% 47-xx: Operator theory
% 47A52: Ill-posed problems, regularization
% 47N20: Applications to differential and integral equations
% 47N40: Applications in numerical analysis


\bigskip
\section{Introduction}
%In \cite{volkov2020stochastic}, we introduced 
%a numerical method for mixed linear and nonlinear inverse problems.
Many physical phenomena are modeled by governing equations 
that depend linearly on some terms and non-linearly on other terms.
For example, the wave equation may depend linearly on a forcing term
and non-linearly on the medium velocity.
%This paper is on inverse problems where both a linear part and a nonlinear
%part are unknown.
 Such inverse problems  occur
in passive radar imaging, or in seismology where the source of an earthquake 
has to be determined (the source could be a point, or a fault) and 
a forcing term supported on that source is also unknown.
The inverse problem is then linear in the unknown forcing term and
nonlinear in the location of the source.\\
In this paper, 
we establish Lipschitz  stability properties for a related class of inverse problems
which we introduce in section \ref{Statement}.
In that class, the associated direct problem is formulated 
by an integral operator $\ca_m$ depending non-linearly on a parameter $m$
and operating on a function $u$. 
In the inversion step, both $u$ and $m$ are unknown but we are only 
interested in recovering $m$. We discuss 
in section \ref{Examples}
examples of such inverse problems
for the  elasticity equation with applications to seismology and for the 
inverse scattering problem in electromagnetic theory.
Assuming a few injectivity and regularity properties  for $\ca_m$,
we prove in section  \ref{Analysis} that the  
 inverse problem 
solved from a  finite number of data points $(\ca_m u (P_{j}))_{1\leq j \leq M}$
has a unique solution  and that this solution is 
Lipschitz stable in the data (theorem  \ref{main disc theorem}).  \\
Since this inverse problem is solvable we can define a function $\psi$
from the vector of data points in $\RR^M$ to $m$.
Under the general assumptions introduced in
section  \ref{Statement},
$\psi $ is Lipschitz continuous (with some insight on its Lipschitz norm
discussed in section \ref{constant estimates}).
Neural networks have been used
as a tool
 for function approximation for some time.
In our case the function $\psi$ is defined implicitly and in applications each evaluation of this function may be  computationally expensive. 
It may  thus be  particularly convenient to pre-compute 
the vector $(\ca_m u (P_{j}))_{1\leq j \leq M}$ for many instances of $u$ and $m$
and use these pre-computed values in a learning algorithm that 
will determine the weights of a neural network ${\cal N}$ approximating
$\psi$.  
It is however important to have some insight about 
how many 
layers and weights 
 ${\cal N}$ should contain. 
Recently,  asymptotic estimates on the size of neural networks 
have been derived
 \cite{yarotsky2017error, shen2021neural, de2021approximation}.
%  \cite{yarotsky2017error requires ReLU or any piecewise linear
% activation function
% corrolary 5.4
 These estimates involve the 
desired accuracy of the approximation, the dimension of the space where
the function is defined, and importantly, the regularity of the function to be approximated. 
In particular, these estimates hold  for Lipschitz regular functions and show that the
accuracy increases in concert 
with the regularity of the function. However, high-order accuracy may be
of little value in practical problems where only noisy inputs are available.
% high order regularity worthless becauseof noise or paucity
In section  \ref{Application} we show an example illustrating how 
a neural network ${\cal N}$ can approximate a  function
$\psi$, where $\psi$ solves a passive  inverse elasticity problem. 
$\psi$  is known to be Lipschitz continuous thanks to 
the theory developed in section  \ref{Analysis}.
This problem 
relates to a model in seismology. Applying ${\cal N}$  to a data
vector $(\ca_m u (P_{j}))_{1\leq j \leq M}$ representing measurements of surface 
displacements yields  $m$ which stands for the geometry parameter
for a planar fault. While preparing data for the learning step and the learning itself
are particularly long and costly (on the order of several hours on a
high performance platform with multiple CPUs), applying  ${\cal N}$
is fast (a few hundredths of a second for 500 evaluations).
 Our computational work shows that
${\cal N}$ performs well on noisy data  too. If more information
on uncertainty is needed
in  the case of noisy data, then the output from ${\cal N}$ 
can be used as a starting point for  
a sampling algorithm aimed at evaluating the covariance of $m$, or possibly its 
probability distribution function. The author proposed 
in  \cite{volkov2020parallel} a parallel sampling
algorithm that
alternates  computing proposals in parallel and combining proposals to accept or reject them.% as in \cite{calderhead2014general}. 
This algorithm, inspired by \cite{calderhead2014general},
 is well-suited to inverse problems 
mixing linear and nonlinear terms, 
where some unknown amount of regularization is necessary, and
where
proposals are expensive to compute. The results 
from \cite{volkov2020parallel}
compare favorably to those obtained from the Maximum
Likelihood (ML), the Generalized Cross Validation (GCV), or the Constrained Least
Squares (CLS) algorithms.

\section{Statement of inverse problem} \label{Statement}
\subsection{Notations and assumptions}\label{assumptions}
Let $R$ and $V$ be two  compact  
manifolds  embedded in $\RR^d$. 
%We assume that 
%\begin{enumerate}[label=(\subscript{M}{{\arabic*}})]
%\item \label{M C2}
%$V$ is of class $C^2$.
%\item \label{M paths}
%$V$ is pathwise connected and there is a constant $L$ such that
%for all $x_1$ and $x_2$ in $V$, there is a path in $V$ connecting $x_1$
%and $x_2$ with length less or equal to $L$.
%\end{enumerate}
%Let 
%\bean \label{generalAm}
%\ca_{\bm} &:& L^2 (R) \ri L^2(V) 
%\eean
Let $\bH_m: V \times R \ri \RR$ be an integration kernel  
depending  on  a parameter $m$, 
where $m$ is
in a compact subset $B$ of $  \RR^p$.
Let   $\ca_{m} $ be the operator defined by convolution by  $\bH_m$
\bean \label{Am1}
 \ca_m g (x) =  \int_R \bH_m(x, y)
 g (y)  d \sigma (y), 
\eean
where $d \sigma (y) $ is the surface measure on $R$.
We assume that $\bH_m$ presents the following regularity properties:

\begin{enumerate}[label=(\subscript{R}{{\arabic*}})]

\item \label{H C1} %\label{Hcont}
$(m,x,y ) \ri \bH_\bm(x, y)$ is continuous  in 
$B \times V\times R $
%\item  \label{kernel}
%All first and second derivatives 
and the gradient in $x$,
$(m,x,y ) \ri \nabla_{x}\bH_m(x, y) 
$, 
%and the Hessian
%and $(m,x,y ) \nabla_{x}\nabla_{x}\bH_m(x, y) 
%$ 
exists and is continuous in
$B \times V \times R $.
\item \label{mder} %(must have $B\subset B'$, $B'$ open for technical reason)
There is an open set $B'$ of $\RR^p$ such that 
$B \subset B'$ and 
the derivatives
$(m,x,y ) \ri \p_{m_i}\bH_m(x, y) $ and
$(m,x,y ) \ri \p_{m_i} \nabla_x\bH_m(x, y) $,
%and  $(m,x,y ) \ri \p_{m_i} \p_{m_j} \bH_m(x, y) $
exist and are continuous in
$B' \times V \times R $, where $1\leq i, j  \leq p$.
\end{enumerate}
Given \eqref{Am1} and assumption \ref{mder}, 
we can define the directional derivative 
$\p_{\bq} \ca_{m} = \nabla_m \ca_m \cdot \bq$ 
of the operator $\ca_{m}$ for any unit vector $\bq$ of $\RR^p$ and $m$ in $B'$. \\

We make the following uniqueness assumptions: 
\begin{enumerate}[label=(\subscript{U}{{\arabic*}})]
\item \label{U1} For any $m, m'$ in $B$ and any $u, v$ in
$H^1_0 (R)$, if  $\ca_{m} u= \ca_{m'} v $ in $L^2(V)$ and
$u \neq 0$ or $v \neq 0$,
then $m=m'$ and $u=v$. In particular, $\ca_m$ is injective for all
$m$ in $B$. 
%\item \label{der prop}  The operator $\nabla_m \ca_{m} \cdot
%\bq$, henceforth simply denoted by $\p_{\bq} \ca_{m}$ is 
%such that for all
% unit vectors $\bq$ of $\RR^p$, all $m$ in $B$, and all $u$ in $H^1_0(R)$,
 %$\p_{\bq} \ca_m u =0 $
%implies that $u =0$.
%injective for all
 %unit vectors $\bq$ of $\RR^p$ and all $m$ in $B$.
\item \label{inter}  
 For all
 unit vectors $\bq$ of $\RR^p$, all $m$ in $B$, and all $u, v$
in $H^1_0(R)$,  $\p_\bq \ca_m  u =  \ca_m v$ implies $u=v=0$. 
In particular $\p_\bq \ca_m$ is injective.
%The intersection of 
%range of $\p_\bq \ca_m $ and of $ \ca_m $is reduced to zero,
%for all $ m $ in   $B$  and all unit vectors $\bq$ of $\RR^p$.
\end{enumerate}
We also assume that integrals over $V$ 
can be approximated by a quadrature of order 1. More precisely,
there is an increasing sequence of integers $M_n$ such that for each $n$ there are  $M_n$ points
$P_j^n, j=1, ..., M_n$, in $V$ and  $M_n$ coefficients 
$C'(j,n), j=1, ..., M_n$, in $\RR$ such that,
\begin{enumerate}[label=(\subscript{Q}{{\arabic*}})]
 \item \label{Q1} For all $\phi$
 in $C^1(V)$,
\bean \label{rule order}
|\sum_{j=1}^{M_n} C'(j,n) \phi (P_{j,n}) - \int_V \phi | =
O(\f{1}{M_n^{\f{1}{\dim V}}}) \sup_{V} |\nabla \phi| .
\eean  
	%where $D^2 \phi$ is the Hessian of $\phi$
	%\item The weights $C'(j,n) $ satisfy 
	%??? really needed ???  or works out instead of and???
	%\bean C'(j,n) > 0
	%\mbox{ and } C'(j,n) = O(\f{1}{M_n}), 
	%\mbox{ uniformly in } j
	%\label{Cjn assum}
	%\eean
	%		\item $\dim H_p  \geq 3N$ and the $3N$ eigenvalues of the discrete operator
%		$\ca_{\bm}^{p, N}(\ca_{\bm}^{p, N })'  $ are equal to 
%		the first $3N$ eigenvalues of the symmetric and compact  operator
%		$\ca_{\bm }\ca_{\bm }'$.
\end{enumerate}
\textbf{Remark:} \\
At first sight, assumption  \ref{inter} may seem unusual. However, useful examples of
operators ${\ca_m}$ satisfying \ref{inter}  abound in inverse problems settings.
Referring to the author's previous work, we can point to an example involving
the Laplace operator in \cite{volkov2021stability}, p 11507 (equation 48 in that paper and subsequent argument), 
and another example involving the elasticity operator
\cite{triki2019stability} p14 (where it is shown that equation 4.5 of that paper imply that $h_0$
and $g_0$ are zero).\\
%We now provide a  more straightforward example which relies on Rellich's lemma 
%and on basic integral operator theory for the Helmholtz operator.
%Set 
%\bea
%\Phi(x,y)=\f{1}{4 \pi} \f{e^{i k|x-y|}}{|x-y|},
%\eea  
%for $x,y$ in $\RR^3$. $\Phi$ is a free-space fundamental solution for the Helmholtz operator. 
%Define the half sphere
%\bea
%R= \{ x \in \RR^3: |x|=1, x_3 \geq 0 \}.
%\eea
%Let $m$ be in the compact interval $B=[0,1]$, and let $V$be the larger sphere centered at the
%origin with radius 3.
We now provide a  more straightforward example which relies on potential theory.
Set 
\bea
\Phi(x,y)=\f{1}{4 \pi} \f{1}{|x-y|},
\eea  
for $x,y$ in $\RR^3$. $\Phi$ is the free-space fundamental solution for the Laplace operator. 
Define the half sphere
\bea
R= \{ x \in \RR^3: |x|=1, x_3 \geq 0 \}.
\eea
Let $m$ be in the compact interval $B=[0,1]$, and let $V$ be the larger sphere centered at the
origin with radius 3. Set $\bH_m=\Phi(x,y+(0,0,m))$ and define 
%\bea %\label{Am1}
 %\ca_m g (x) =  \int_R \bH_m(x, y)
 %g (y)  d \sigma (y), 
%\eea
$\ca_m $ by \eqref{Am1} for $g$ in $H^1_0(R)$ and $x$ in $V$.
Now assume that
$\p_\bq \ca_m  u =  \ca_m v$ on $V$ for some $m$ in $B$ and $u, v $ in $H^1_0(R)$,
and some $\bq$ such that $|\bq | =1$.
Note that 
\bea
\p_\bq \ca_m  u =
\epsilon \int_R \p_{y_3}\Phi(x,y+(0,0,m))
 u (y)  d \sigma (y),
\eea
where $\epsilon = \pm 1$.
Let $R_m = R+(0,0,m)$. For  $x$ in $\RR^3 \sm R_m$, introduce the function,
\bea
w (x)=  
 \epsilon \int_R \p_{y_3}\Phi(x,y+(0,0,m))
 u (y)  d \sigma (y)
- \int_R \Phi(x,y+(0,0,m))
 v (y)  d \sigma (y),
\eea
  As $\Delta w = 0$ in $\RR^3 \sm R_m$, $w$ 
decays at infinity, 
and $w=0$ on the sphere $V$ (this is due to the assumption $\p_\bq \ca_m  u =  \ca_m v$), %as $w$ decays at infinity,
since the exterior Dirichlet problem has at most one solution, we claim
 that $w$ is zero in the exterior of the sphere $V$.
Since $w$ is analytic in  $\RR^3 \sm R_m$, it must equal zero everywhere 
in $\RR^3 \sm R_m$. We now use the spherical coordinates $(r, \theta, \phi)$, where $0\leq \theta
\leq \pi$ is the co-latitude. Then for  $y$ on the half sphere $R$, 
\bea \p_{y_3} = \cos \theta \f{\p}{\p n}  - \sin \theta \p_\theta,
\eea
where $n$ is the exterior normal vector on $R$. Note that $\p_\theta$ is a tangential
derivative on $R$.
% and $\omega(\theta, \phi)$ is a smooth tangential field.
Now, applying well-known properties of the single and double layer potential,
we find that 
\bea
\lim_{h\ri 0, h >0} w(x+ h n) - w(x-hn) = \epsilon \cos \theta u(x),
\eea
where $x$ on the half sphere $R_m$, and $n$ is the exterior normal at $x$.
As $w$ is zero  in $\RR^3 \sm R_m$, it follows that $u$ is zero.
We can  now show  likewise that $v$ is zero by using the jump formula 
for the normal derivative of the single layer potential.

\subsection{Statement of the continuous and 
the discrete inverse problems} \label{prob statement}
The continuous inverse problem:\\
\begin{center}
Given $\ca_{m} u $ in $L^2(V)$ for some unknown $u$ in $H^1_0(R)$
and $m$ in $B$, find $m$.
\end{center}%\\\\
The discrete  inverse problem:\\
\begin{center}
Given $n$ in $\NN$ 
and the discrete values
$\ca_{m} u (P_{j,n})$, $j=1,.., M_n$, for some unknown $u$ in $H^1_0(R)$
and $m$ in $B$, find $m$. 
\end{center}%\\\\

It is clear from \ref{U1} that the continuous inverse problem has a unique solution.
Combining some of the other assumptions, we will show that  the continuous inverse problem 
is in some sense Lipschitz stable. 
Interestingly, \ref{U1}  also implies that there is a unique $u$ producing $\ca_{m} u $.
However, since the linear operator $\ca_{m} $ is compact, 
 $u$ does not depend continuously on  $\ca_{m} u $. 
Solving the linear inverse problem consisting of estimating $u$  from $\ca_{m} u $
is a classical problem once the nonlinear parameter $m$ is known and will not
be covered in this paper.
We will  also show that for all $n$ large enough the discrete inverse 
problem is  uniquely solvable and Lipschitz stable as well.

\section{Examples of inverse problems satisfying the conditions 
stated in section \ref{assumptions}
} \label{Examples}
% Laplace, elasticity, old papers by Colton / Kress
\subsection{Fracture or wall in half space governed by the Laplace equation}
\label{fracture}
%Here, we examine  the case of a model involving the Laplace equation.
	This model %involving the Laplace equation 
	is  relevant to geophysics: in  dimension two,
	it relates to the so called anti-plane strain configuration.
	This  configuration has  attracted much attention 
	from geophysicists and mathematicians due to 
	how simple and yet relevant
	  this formulation is
	\cite{dascalu2000fault, ionescu2006inverse, ionescu2008earth}.
 In dimension three, this model  relates to irrotational incompressible 
flows in a medium with a top wall and an inner wall.\\
Let $\RR^{3-}$ be the open half space  $\{ x_3<0 \}$, where
we  use the notation $x=(x_1,x_2,x_3)$ for points in $\RR^3$.
  Let $\Gamma$ be a Lipschitz  open surface in $\RR^{3-}$ and $D$ a domain in $\RR^{3-}$
	with Lipschitz boundary such that $\Gamma \subset \p D$.
	We assume that $\Gamma$ is strictly included in $\RR^{3-}$ so that  the distance 
	from $\Gamma$ to the plane $\{ x_3 =0 \}$ is positive.
We define the direct fracture (or crack) problem to be the boundary value problem, 
	\bean
        \Delta u=0\text{ in }\RR^{3-}\setminus \ov{\Gamma},  \label{BVP1}     \\
        \p_{x_3} u=0\text{ on  the surface } x_3=0,  \label{BVP2}     \\
     \aaa \frac{\p u}{\p \bn} \bbb = 0 \mbox{ across }\Gamma,   \label{BVP3}\\
      \aaa u \bbb
			=g \mbox{ across }\Gamma,  \label{BVP4}\\
    u(x)=O ({\frac{1}{| x|}}) \text{ uniformly as } |x|\to \infty,     \label{Decay1}
\eean
where  $\aaa v \bbb $ denotes the jump of a function $v$ across $\Gamma$ in the normal
direction,
and $\bn$ is a unit normal vector  to $\Gamma$.
In problem (\ref{BVP1}-\ref{Decay1}),
$u$ can be a model for the potential of an irrotational flow,
with an impermeable and immobile wall $\{x_3  =0 \}$. 
%The flow is tangent to the inner wall
%$\Gamma$, and 
The discontinuity of the tangent flow is given by the 
tangential gradient of $g$. Alternatively, if we wrote the analog of problem
 (\ref{BVP1}-\ref{Decay1}) in two dimensions, it could model a strike-slip fault in geophysics
where the equations of linear elasticity simplify to the scalar Laplacian.
In that case the scalar function $u$
models displacements in the direction orthogonal to a cross-section,
 $g$ models a slip, and normal derivatives model traction
\cite{dascalu2000fault, ionescu2006inverse, ionescu2008earth}. 
Mathematically, problem (\ref{BVP1}-\ref{Decay1}) can be stated 
for $u$ in the functional space 
\bea
\{ v \in H^1_{loc} (\RR^{3-}\setminus \ov{\Gamma}): \nabla v,
 \f{v}{\s{ 1 + |x|^2}} \in
   L^2(\RR^{3-}\setminus \ov{\Gamma})\},
\eea
and
uniqueness 
%for problem (\ref{BVP1}-\ref{Decay1})
can be shown by setting up a variational problem for $u$,
see \cite{volkov2021stability}.
For existence, we can choose 
to express $u$ as an integral over $\Gamma$ against an adequate Green function. %since this formulation will be crucial
%further in this study. 
Indeed, denoting
\bean \label{fund}
\Phi(x,y) = \f{1}{4 \pi} \f{1}{|x - y|},
\eean
 the free space Green function
for the Laplacian, we can argue that
\bean \label{freeint}
 u'(x) = \int_{\Gamma} \nabla_{y} \Phi (x,y)  \cdot \bn(y) g(y) d \sigma (y),
\eean
satisfies \eqref{BVP1}, \eqref{BVP3}, \eqref{BVP4}, \eqref{Decay1}.
%, and is in ${\cal V}$.
%Indeed, as $g$ is in $ \tilde{H}^{1/2}(\Gamma)$, we can extend $g$ by zero to a function
%in $ H^{1/2}(\p D)$ so (\ref{freeint}) equals
%\bean \label{freeint2}
% u'(x) =  \int_{\p D} \nabla_{y} \Phi (x,y)  \cdot \bn(y) g(y) d \sigma (y).
%\eean
Even if  $\p D$ is only Lipschitz regular, by Theorem 1 in \cite{costabel1988boundary},
%(\ref{freeint2}) defines 
$u'$ is a function in $H^1_{loc} (\RR^{3} \setminus \ov{\Gamma})$
and by Lemma 4.1 in \cite{costabel1988boundary} the jump $[ u']$ across
$\Gamma$ is equal to $g$ almost everywhere, while the jump $[ \f{\p u' }{\p \bn}]$ is zero.
To find a solution to the PDE (\ref{BVP1}-\ref{BVP4})
we then set
\bean \label{intform}
 u(x) =  \int_{\Gamma} \bH (x,y,\bn)   g(y) d \sigma (y),
\eean
where 
\bean \label{Hformula}
\bH (x,y,\bn) = \nabla_{y} \Phi (x,y)  \cdot \bn(y)  +
\nabla_{y} \Phi (\ov{x},y)  \cdot \bn(y) ,
\eean
and $\ov{x} = (x_1, x_2, - x_3)$.
Then conditions (\ref{BVP1}-\ref{BVP4}) are satisfied,  $u(x) = O(\f{1}{|x|^2})$,
 and $\nabla u(x) = O(\f{1}{|x|^3})$, uniformly in 
$\f{x}{|x|}$. %, so $u$ is in ${\cal V}$. 
In \cite{volkov2021stability}
we considered the case where $\Gamma$ is planar and in that case
$m$ is a geometry parameter in $\RR^3$.
In fact, $\Gamma$ could  also be included in several contiguous polygons
and as a result $m$ would be in a higher dimensional space. 
There is a limitation on the number of polygons as suggested 
by the counterexample provided in appendix A of 
\cite{volkov2021stability}.\\
For the planar case, the dependance of $\Gamma$ on the geometry parameter
$m$ can be modeled as follows. 
Let $R $ be the closure of a   bounded relatively open  set 
 in the plane $\{x_3=0 \}$  and
$m= (m_1,m_2,m_3)$ in $\RR^3$.
Define 
the surface in $\RR^3$
 \bean \label{Gammam def}
\Gamma_{m}=\{(x_1,x_2,m_1 x_1+m_2x_2+m_3 ):(x_1,x_2)\in R \}.
\eean 
 Let $B$ be a closed and bounded set of $m$ 
in $\RR^3$
such that $\Gamma_{m}\subset\RR^{3-}$.
There is a negative constant $d_0$ 
such that 
%We assume that $B$ is  . 
%It follows  that
%there is a positive constant $\beta$ such that
\begin{equation}
    \label{BProperty}
  m_1 x_1+m_2x_2+m_3\leq d_0, \q \forall (x_1,x_2)\in R.
\end{equation}
%for all $\bm$ in $B$, where $\mbox{dist }$ is the distance between sets.
We choose the unit normal vector on $\Gamma_m$ to be
%$\bn=\frac{(-a,-b,1)}{\sqrt{a^2+b^2+1}}$
 $\bn=\frac{(-m_1,-m_2,1)}{\sqrt{m_1^2+m_2^2+1}}$. The surface element on
$\Gamma_m$ will be
 denoted by
%$\sigma d x_1 d x_2=\sqrt{a^2+b^2+1} \, d x_1 d x_2$. 
$\sigma d x_1 d x_2=\sqrt{m_1^2+m_2^2+1} \, d x_1 d x_2$. 
%Let $H_0^1(R)$ be the space of functions $g$ on $R$ with Sobolev $H_0^1$ regularity. Let $V$ be a relatively open subset of $\{x_3=0\}$. 
 According to \eqref{intform} and \eqref{Hformula}, we can define the 
operator
%introduce the linear compact operator mapping jumps $g$ across 
%$\Gamma_m $ to Dirichlet data on $V$,
\bean
 \ca_m:H_0^1(R)\to L^2(V), \no \\
      (\ca_m g )(x) =\int_R \bH(x_1,x_2,0,y_1,y_2,m _1 y_1+m _2 y_2+ m_3,\bn)
				g(y_1,y_2)\sigma
				dy_1dy_2,    \label{AmDef}
\eean
where $V$ is the closure of a bounded open set of the top plane
 $\{x_3=0\}$.
Due to \eqref{Hformula} and \eqref{BProperty}
 $\bH$ is  analytic in $(m,x,y) \in B \times V \times R$,
 thus $\bH$ satisfies \ref{H C1} and
\ref{mder}.
We know from Theorem 2.2 in \cite{volkov2021stability}  
  %\ref{InverseProblemResult} 
that $\ca_m$ is injective.
 %Fix a  non-zero $h$ in $H_0^1(R)$ and define the 
%function $\phi:B\to L^2(V)$ by 
%\begin{equation} 
  %  \label{phiDef}
   % \phi(\bm)= A_m h.
%\end{equation}
%Theorem \ref{InverseProblemResult} implies that $\phi$ is injective. 
%We  showed in
%\cite{volkov2021stability}
 %that the inverse of $\phi$ defined on $\phi(B)$ and valued in $B$ is of class $C^1$ by applying
%the  inverse function Theorem.
%As a consequence, $\phi^{-1}$ is Lipschitz continuous, which will yield
%our first stability estimate.
In addition,
 assumption \ref{U1}
holds due to the same theorem.
The proof of theorem 2.3 in \cite{volkov2021stability}
(in particular equation (26) from that paper leading to $\gamma_1=
\gamma_2=\gamma_3=0$) shows that
$\p_{\bq} \ca_m  $ is injective for all $m$ in $B$ and $\bq$ in $\RR^3$
with $|\bq|=1$. 
Assumption \ref{inter}  is also satisfied: this is shown in page 11507
 of \cite{volkov2021stability}, right beneath equation 48 of that paper.

\subsection{Fault in elastic half space}\label{sec Fault in elastic}
%Using the standard rectangular coordinates $x= (x_1, x_2, x_3)$ of $\RR^{3}$,
%A point $(x,y,z)$ will denote a point
%we define $\RR^{3-}$ to be the open half space $x_3<0$.
%Let  $\Gamma$ be a bounded open surface
%with smooth boundary included in the plane $\{ x_3=0 \}$.
%The derivative in the $i$-th coordinate will be denoted by $\p_i$.
Assume that $\RR^{3-}$ is a linear isotropic elastic medium with 
 Lam\'e constants $\lambda$ and $\mu$ 
such that $\mu >0$, $\lambda + \mu >0$.
Denote $e_1, e_2, e_3$ the natural basis of $\RR^3$.
%will be two positive constants.
For a vector field $u = (u_1, u_2, u_3)$,
%the stress  and strain tensors will be denoted as follows,
%\bea
%\sigma_{ij}(u) = \lambda \, \di u \, \delta_{ij} + \mu \, (\p_i u_j + \p_j u_i ), \\
%\epsilon_{ij}(u) = \f12 (\p_i u_j + \p_j u_i ),
%\eea
denote the   stress vector in the normal
direction $\bn$, 
$$\ds T_n u = \sum_{1 \leq j \leq 3} (\lambda \, \di u \, \delta_{ij} + \mu \, (\p_i u_j + \p_j u_i )) n_j e_i.$$
Let $\Gamma$ be a Lipschitz open surface which is strictly included in $\RR^{3-}$.
Let $u$ be the displacement field solving
\bean
\mu \Delta u+ (\lambda+\mu) \nabla \di u= 0  \mbox{ in } \RR^{3-}
\setminus \Gamma \label{uj1}, \\
 T_{e_3} u =0 \mbox{ on the surface } x_3=0 \label{uj2}, \\
 T_{\bn} u  \mbox{ is continuous across } \Gamma \label{uj3}, \\
 \aaa  u\bbb =g \mbox{ is a given jump across } \Gamma , \label{uj3andhalf}\\
u (x) = O(\f{1}{|x|^2}),  \nabla u (x) = O(\f{1}{|x|^3}), \mbox{ uniformly as }
|x| \rightarrow \infty,
\label{uj4}
\eean
%where %$x=(x_1, x_2, x_3)$ is any point in $\RR^{3-}$ and
%$e_3$ is the vector $(0,0,1)$.
In \cite{volkov2017reconstruction}, we proved that
 problem
(\ref{uj1}-\ref{uj3andhalf})
has a unique solution in  the functional space of vector fields $u$
such that  $\nabla u $ and in $\f{|u|}{ (1 + |x|^2 )^{\f12}}$ are 
in $L^2(\RR^{3-} \setminus \ov{\Gamma})$.
In particular, this applies to the case where $g$ is in $H^1_0(\Gamma)$.
%Theorems  \ref{direct exist unique} and \ref{uniq1} 
%were  proved in \cite{volkov2017reconstruction} 
%for media with constant Lam\'e coefficients. Later,
In \cite{aspri2020analysis}, 
the direct problem (\ref{uj1}-\ref{uj3andhalf})
 was analyzed under weaker regularity conditions 
for $u$ and $g$. In \cite{aspri2020arxiv},
the direct problem (\ref{uj1}-\ref{uj3andhalf})
was proved to be uniquely solvable in case of
piecewise Lipschitz coefficients and general elasticity tensors.
Both \cite{aspri2020arxiv}  and \cite{aspri2020analysis}
include a proof of  uniqueness  for the  fault inverse problem under appropriate assumptions.
%In addition, the solution $u$ satisfies the decay conditions
%(\ref{uj4}).
%\end{thm}
%In this paper we will only consider forcing terms $g$ which are tangential to
%$\Gamma$. Physically, this reflects that the fault $\Gamma$ is not opening or starting to self intersect: only slip is allowed.
 %We recall that if $g$ is continuous, the support of $g$,
%$\mbox{supp } g$, is equal to the closure
%of the set of points in $\Gamma$ where $g$ is non zero; in general $\mbox{supp }  g$
%is defined in the sense of distributions.
However, we need to use in this paper the more regular framework 
introduced in \cite{volkov2017reconstruction} since we have to use
the stability properties derived in \cite{triki2019stability} and they require
$g$ to be in $H^1_0(\Gamma)$. \\
There is a Green's tensor  $\bH$ such that if $g$ is in $H^1_0(\Gamma)$,
 		the  solution $u$ to problem (\ref{uj1}-\ref{uj3andhalf})  can also be written out as the convolution on $\Gamma$
\bean
  \int_\Gamma \bH(x, y) g(y) \, d \sigma (y) \label{int formula}.
\eean
The practical determination of this adequate half space Green's tensor $\bH$ was  studied
in  \cite{Okada} and later, more rigorously, in \cite{volkov2009double}.
Suppose that $\Gamma$  is such that it can be parametrized
 by  $m$ in $\RR^p$.
In  \cite{volkov2020parallel}, $\Gamma$ was modeled as two contiguous
quadrilaterals with known first two coordinates and accordingly $m \in \RR^6$.
In  \cite{volkov2020stochastic}, $\Gamma$ was modeled as a
parallelogram which is the projection of a rectangle on
an unknown plane so $m \in \RR^3$: accordingly, 
 $\Gamma_m$ can be defined as previously by \eqref{Gammam def}.
For sake of simplicity, assume that this is the case.
We still assume that the distance condition \eqref{BProperty}
is satisfied. 
We thus obtain displacement vectors for $x$ in $V$ by the integral formula
\bean \label{new int}
u(x) = \int_R \bH_m(x, y_1, y_2)
 g (y_1, y_2) \sigma d y_1 d y_2,
\eean
for any $g$ in $H^1_0 (R)$ and  $m$ in $B$, where $\sigma $ is the surface element on $\Gamma_m$
and  $\bH_m(x, y_1, y_2)$ is derived from
the Green's tensor $\bH$ for $y$ on $\Gamma_m$.
%Let $V$ be a non-empty, open subset of the plane
%$x_3 =0 $ and
%We now assume that $V$ is a bounded open subset of the plane $x_3 =0$ and for
%a fixed $\tilde{u}$ be in $L^2(V)$,
%and  a fixed $m$ in $B$ we define the functional
Define the operator
\bean
\ca_{m} &:& H^1_0 (R) \ri L^2(V) \no \\
&&g \ri \int_R  \bH_m(x, y_1, y_2)
 g (y_1, y_2) \sigma d y_1 dy_2 .  \label{Aabd elastic}
\eean
Note that in this model both $g$ and $\ca_m g$ are vector fields.
 Anyway, the generalization to vector fields of the assumptions made in section
 \ref{assumptions} for scalar functions is straightforward. 
%It is clear that $A_{m} $ is linear, continuous, and compact.
The closed formula for $\bH (x,y)$ is involved \cite{volkov2009double} 
but it is a  real analytic function of $(x,y)$ if $x_3 \leq 0$, $y_3<0$, 
and  $x  \neq y$. Consequently, thanks to
the distance condition \eqref{BProperty}, assumptions 
\ref{H C1} and  \ref{mder} are satisfied. 
Uniqueness assumption \ref{U1} was proved in 
\cite{volkov2017reconstruction}, theorem 2.1.
The proof of theorem 3.1 in \cite{triki2019stability}
(in particular equation (3.20) from that paper leading to $\gamma_1=
\gamma_2=\gamma_3=0$) shows that
$\p_{\bq} \ca_m  $ is injective for all $m$ in $B$ and $\bq$ in $\RR^3$
with $|\bq|=1$. The argument holds if the slip $g$ is tangential 
on $\Gamma_m$. 
Similarly, assumption \ref{inter} holds thanks to a result shown in 
\cite{triki2019stability}. Starting from equation (4.5) in \cite{triki2019stability},
it is shown that 
if $\p_{\bq} \ca_m h = \ca_m g$ for some $g$ in $H^1_0(R)$, then
$h$ must be zero. This was done 
under the additional assumption that $h$
% (or $g$ in equation 
%\eqref{uj3andhalf}) 
is either one-directional or the gradient of a function
with Sobolev regularity $H^2$, while still in $H^1_0$.\\
%It was assumed in  \cite{triki2019stability}
%that the set $B$ is such that for $m \in B$, $\Gamma_m$ is in the half-space
%$\RR^{3-}$. 
%This assumption may prove inadequate in practical calculations.
%Instead, it may be more convenient to fix a depth $d_0 <0$ 
%and to only take into account points $(x_1, x_2, a x_1 + b x_2  + d )$
%for $(x_1, x_2)$ in $R$ such that 
%$a x_1 + b x_2  + d  \leq d_0$  in the following manner:
%and define
%$\Gamma_m$ to be
%$$
%\{ (x_1, x_2, a x_1 + b x_2  + d ): (x_1, x_2) \in R \mbox{ and } 
%a x_1 + b x_2  + d  \leq d_0\}.
%$$
%If $R$ is a rectangle
%$g$ defined on subspace .... further study!
%but numerics are here
%fix $\phi$ in $C^{\infty} (\RR)$ such that $0\leq \phi \leq 1$,
%$\phi(t) = 1$, if $t \leq 2 d_0$,  
%$\phi(t) = 0$, if $t \geq d_0$.
%Define
 %the operator
%\bea
%\ca_{m} &:& H^1_0 (R) \ri L^2(V) \no \\
%&&g \ri \int_R  \bH_m(x, y_1, y_2)
%\phi( a y_1 + b y_2 +d) g (y_1, y_2) \sigma d y_1 dy_2 . % \label{Aabd elastic}
%\eea
%This time, if $\ca_m g =0$, then $\phi g =0$.
 %It is clear that  assumptions \ref{H C1}, \ref{mder}, \ref{U1}.
%still hold
%for this more sophisticated model.
%However, proving that assumption \ref{inter}  holds requires more work:
%we plan to do that in future work.
% we will revisit 
%the calculations and arguments from \cite{triki2019stability}
% could be a great problem for Yulong

\subsection{Inverse acoustic scattering problem}%{Colton Kress waves}
Let $D$ be a Lipschitz domain in $\RR^3$ modeling a soft scatterer
for acoustic waves. 
When this scatterer is illuminated by a plane wave $e^{i k \omega \cdot x }$,
where $k>0$ is the wavenumber,
it produces a scattered field $u$ which satisfies the following PDE:
\bean
        (\Delta + k^2) u=0\text{ in }\RR^{3-}\setminus \ov{D},  \label{w1}     \\
        u(x)=- e^{i k \omega \cdot x } \text{ on  the surface }\p D,  \label{w2}     \\
         \nabla u \cdot \f{x}{|x|}  -iku =O ({\frac{1}{| x|^2}}) \text{ uniformly as } |\bx|\to \infty.     \label{Decay w}
\eean
It is well known that  problem (\ref{w1}-\ref{Decay w}) is uniquely solvable
and that the solution $u$ satisfies as $| x | \ri \infty$,
\bea
u(x) = \f{e^{ik |x|}}{|x|} ( u_\infty (\hat{x})+ O(\f{1}{|x|^2}) ) ,
\eea
where $\hat{x} = \f{x}{|x|}$ is a variable on the unit sphere 
$S$ of $\RR^3$. $u_\infty$ is called the far-field pattern of $u$. \\
The inverse scattering problem consists of reconstructing $D$ from the 
far-field pattern $u_\infty$. 
With additional assumptions,
it is known in a few cases that  $D$ can be determined from $u_\infty$. 
For example, Alessandrini and Rondi  proved
that if it is initially known that $D$ is polyhedron,
%this inverse scattering problem has a unique solution
this determination is possible
\cite{alessandrini2005determining}.
For general shapes, it was proved 
that if for infinitely many directions $\omega$ of 
of the
incident plane wave
the far field $u_{\infty, \omega} $ is given, then $D$ is uniquely
determined (\cite{colton1998inverse}, theorem 5.1).
Using infinitely many incident plane waves 
may be prohibitive in practice, but
interestingly, it was shown  if $D$ is included in a ball of radius $r$
and  the wavenumber $k$ satifies $kr < \pi$, 
$D$ can be again determined from $u_\infty$
(\cite{colton1998inverse}, corollary 5.3).\\
Note that in this classical inverse scattering problem, 
the forcing term $e^{i k \omega \cdot x} $ in 
\eqref{w2} is a known incoming wave. If we know 
that the geometry of $D$ can be parametrized by some 
$m$ in $\RR^p$ (for example if $D$ is known to be 
a polyhedron or an ellipsoid)
the classical inverse scattering problem
is therefore 
not as general as the problem stated in 
\eqref{prob statement} where $u$ is unknown.
A closely related model that uses the full generality 
of  the problem stated in 
\eqref{prob statement} corresponds to the case where an unknown wave 
illuminates $D$. This unknown wave may not necessarily be a plane wave.
 Proving that $B$ can be determined from 
$u_\infty$ 
 in this more challenging case too will be the subject of future work.\\
If $m$ is a parameter determining the geometry of $D$ and 
$\ca_m$ is the operator mapping the incoming wave $e^{ik \omega \cdot x}$
to the far-field $u_\infty$, it is known that $\ca_m$ is injective 
(if $D$ is a polyhedron or an ellipsoid while $k$ is small enough).
The differentiability $\ca_m$ in $m$ 
and the injectivity of $\p_\bq \ca_m$
have been established,
see theorems 5.14 and 5.15 in \cite{colton1998inverse} 
and  \cite{potthast1994frechet}.
 % Potthast [273, 274].
%273. Potthast, R.: Fr´echet differentiability of boundary integral operators in inverse acoustic scattering.
%Inverse Problems 10, 431–447 (1994).
%274. Potthast, R.: Fr´echet Differenzierbarkeit von Randintegraloperatoren und Randwertproblemen
%zur Helmholtzgleichung und den zeitharmonischen Maxwellgleichungen.%
%Dissertation, G¨ottingen 1994.
 However, the full scope of assumption  \ref{inter}  will have to be 
studied in future work.

\section{Analysis and proof of stability results} \label{Analysis}
\subsection{The continuous case}\label{The continuous case}
\begin{lem}\label{intermediate lemma}
Assume that $u_k$ converges weakly to $u$ in $H^1_0(R)$.
Fix $m$ in $B$.
Then $\ca_m u_k - \ca_m u $ converges uniformly to zero in $V$.
Let $m_k$ be a sequence in $B$ converging to $m$. Then 
$\ca_{m_k} u_k- \ca_m u $ converges uniformly to zero in $V$.
% and
%\bean \label{disc conv}
%\sum_{j=1}^{M_{n_k}} C'(j,n_k) | 
%(\ca_{\bm_k}{u_k} - \ca_{\bm}{u} ) (P_j)|^2
%\eean
%converges to zero.
\end{lem}
\textbf{Proof:}
According to  \eqref{Am1} and \ref{H C1},
\bean
|\ca_m u_k (x)- \ca_m u(x)|
&=&|\int_R  \bH_m(x, y)
 (u_k (y) - u (y)) d\sigma (y) | \no \\
&\leq& \sup_{ x \in V, y \in R } | \bH_m(x, y) ||R|^{\f12}
(\int_R  
 (u_k (y) - u (y))^2 d\sigma (y) )^{\f12},
\label{uni estimate}
\eean
and since $u_k$ converges strongly to $u$ in $L^2(R)$,
the first claim is proved.
To prove the second claim, it suffices to show that 
$\ca_{m_k}  u_k- \ca_m u_k$ converges uniformly to zero.
This is due to \ref{H C1} and the estimate,
\bean
|\ca_{m_k} u_k (x)- \ca_m u_k(x)| \leq  
\sup_{ x \in V, y \in R } | \bH_{m_k}(x, y ) -
\bH_m(x, y) ||R|^{\f12}
(\int_R  
 u_k (y) ^2 d \sigma(y) )^{\f12}.
\label{uni estimate 2}
\eean
%and the lemma is proved. \\
%Finally, to prove \eqref{disc conv}
$\Box$\\

We introduce the following notations: for $\phi$ in $H^1_0(R)$, we set
\bean \label{norm R}
\|\phi \| = \big(\int_R |\nabla \phi(y)|^2  d\sigma(y)\big)^{\f12}. 
\eean
For a function $\psi$  in $L^2(V)$,
\bean \label{norm V}
\|\psi \| = \big(\int_V | \psi(x)|^2  d\sigma'(x)\big)^{\f12},
\eean
where $\sigma' $ is the surface element in $V$.\\
%Let $M_1$ and $A_2$ be two positive constants.
%	Define
%		the subset of $H^1_0(R)$, (or THERE IS?)
%		\bean \label{setS}
%	S = \{ \varphi \in H^1_0(R): 
%	\| \varphi \| \leq A_2 \mbox{ and }
%	M_1 \leq \| \ca_{m} \varphi \| 
%	\mbox{ for all } m \in B
%	\},
%\eean

We endow $L^2(V)$ with its usual Hilbert space inner  product structure associated to the norm \eqref{norm V}.
In $H^1_0(R)$, thanks to  Poincare's inequality, 
\eqref{norm R} defines an equivalent norm
associated to the inner product 
$<\phi, \psi> = \int \nabla \phi \cdot \nabla \psi$.
We choose to identify $H^1_0(R)$ with its dual through this inner product. 
It is clear that $\ca_m$ defines a compact linear map from 
$H^1_0(R)$ to $L^2(V)$. Let $\ca_m^*$ be its dual, continuously mapping 
$L^2(V)$ to $H^1_0(R)$. $\ca_m^* \ca_m$ is then a compact and symmetric map 
from $H^1_0(R)$ to $H^1_0(R)$.
By continuity of $\ca_m$ in $m$ and compactness of
 $B$, the minimum of $\| \ca_m^* \ca_m \|$ 
for $m$ in $B$ is achieved. By  \ref{U1}, $\ca_m$ is injective for all 
$m$ in $B$, so $\| \ca_m^* \ca_m \|$  is  in particular  non zero. 
We now fix $\beta >0$ such that $\| \ca_m^* \ca_m \| > \beta^2 $ for all $m$ in $B$.
Let $E_m$ be the subspace spanned by the eigenvectors of $\ca_m^* \ca_m $ corresponding to eigenvalues strictly greater than $\beta^2$. Necessarily,
\bean \label{lower bound}
\forall m \in B, \, \forall v  \in E_m, \q \| \ca_m  v\| \geq  \beta\| v\|. 
\eean
%A compactness argument can show that there is a $\beta' >0$ such that,
%\bean \label{lower bound2} \forall \bq \in \RR^p
 %\mbox{ with } \bq=1, \forall m \in B, \, \forall v  \in E_m, \q \| \p_\bq \ca_m  v\| > \beta\| v\|. 
%\eean
Let $\cp_m $ be the orthogonal projection in $H^1_0(R)$ on 
the finite dimensional space $E_m$. 
%Note that since $E_m$
%is finite-dimensional,  $\cp_m $  is well defined.
\begin{lem}\label{proj lemma}
Let $m_0$ be in $B$. If $\beta^2$ is not an eigenvalue of  
$\ca_{m_0}^* \ca_{m_0} $, the estimate 
\bean \label{proj_est}
\| \cp_m - \cp_{m_0} \| = 0(|m - m_0|)
\eean
 holds. 
\end{lem}
\textbf{Proof:}
Let 
$$\lambda_1^2 \geq ...  \geq \lambda_p^2>\beta^2 >
\lambda_{p+1}^2  \geq  \lambda_{p+2}^2 ...$$
be the eigenvalues of $\ca_{m_0}^* \ca_{m_0} $. 
%Fix $t$ in $(\beta^2, \lambda_r^2)$. 
Let ${\cal C}_1$ be the  circle in the complex plane centered at the origin with radius
 $\lambda_1^2+1$ and ${\cal C}_2$ be the  circle  centered at the origin with radius
 $\beta^2$.
Then $\cp_{m_0}$ can be written as the combination of contour integrals \cite{kato2013perturbation}
\bean \label{int for proj}
\cp_{m_0} = \f{1}{2i \pi} \int_{{\cal C}_1} (zI -  \ca_{m_0}^* \ca_{m_0} )^{-1} dz-
 \f{1}{2i \pi} \int_{{\cal C}_2} (zI -  \ca_{m_0}^* \ca_{m_0} )^{-1} dz. 
\eean
For all $z$ in ${\cal C}_2$, $\| (zI -  \ca_{m_0}^* \ca_{m_0} )^{-1} \|$ is  bounded
by $\max \{ (\lambda_p^2 - \beta^2)^{-1}, (\beta^2 - \lambda_{p+1}^2 )^{-1}\}$
and it is clear that 
for $z$ in ${\cal C}_1$, $\| (zI -  \ca_{m_0}^* \ca_{m_0} )^{-1} \|$ is uniformly bounded.
It follows that 
for $m$ in an open neighborhood of $m_0$,
$ (zI -  \ca_{m}^* \ca_{m} )^{-1}$ is defined and uniformly bounded
for all $z$ in ${\cal C}_1$ and in ${\cal C}_2$. 
The estimate \eqref{proj_est} now results from the factorization
\bean
&& (zI -  \ca_{m}^* \ca_{m} )^{-1}  -  (zI -  \ca_{m_0}^* \ca_{m_0} )^{-1}
\no  \\
&=&
(zI -  \ca_{m}^* \ca_{m} )^{-1}  (\ca_{m}^* \ca_{m}  -
\ca_{m_0}^* \ca_{m_0}
) (zI -  \ca_{m_0}^* \ca_{m_0} )^{-1}, \label{factorization}
\eean
and assumption \ref{mder}.
$\Box$\\

%Fix two positive constants $M_1$ and $A_2$ and define the subset 
%of $H^1_0(R)$
%\bea
%S = \{ \phi \in H^1_0(R) :   \| \phi \| \leq M_1 \mbox{ and }
%\forall m \in B, \| \ca_m \phi \| \geq A_2
%\}.
%\eea
%this contradicts \ref{U1}.\\

\begin{thm}\label{unif with Em}
%Fix two positive constants $A_1$ and $A_2$.
Fix $\beta>0$ and define $E_m $ as previously.
There is a positive constant $C$ such that for all $m, m'$ in $B$, all $u$ in $E_m$ and all $v$ in $E_{m'}$,
\bean \label{main est uni}
	 \| \ca_m u -  \ca_{m'} v  \| \geq
 C  \|\ca_{m'} v\|  |m - m'|.
\eean
\end{thm}
\textbf{Proof:}
Since $\ca_m, \, \ca_{m'}$ are linear operators and $E_m, \, E_{m'}$ 
are linear spaces we only need to show this estimate in the case where
$\|\ca_{m'} v\| =1$. 
Arguing by contradiction, assume that there are two sequences
 $m_k$  and $m_k'$ in $B$ with $m_k\neq m_k'$ for all $k$,   a sequence
$u_k$ in $E_{m_k}$, and a sequence $v_k$ in in $E_{m_k'}$
with $  \| \ca_{m_k'} v_k \|=1\textit{} $
 such that 
	\bean \label{contra_cont3 uni}
	 \| \ca_{m_k} u_k -  \ca_{m_k'} v_k \| < 
 \f{1}{k}  |m_k - m_k'|.
\eean
Given relation \eqref{lower bound}
$v_k$ is  bounded,
 so by  \eqref{contra_cont3 uni} and \eqref{lower bound},
$u_k$ is bounded too.
Without loss of generality, we may assume that $m_k$ converges to some $m$ in $B$, $m_k'$ converges to some $m'$ in $B$,
$u_k$ is weakly convergent to some $u$ in $H^1_0(R)$,
$v_k $ is weakly convergent to some $v$ in $H^1_0(R)$.
By lemma \ref{intermediate lemma}, 
$\ca_{m_k} u_k$ converges strongly to $\ca_{m} u$ and
$\ca_{m_k'} v_k$ converges strongly to $\ca_{m'} v$.
Since $\| \ca_{m'} v \| =1 $,
 combining \ref{U1} and
\eqref{contra_cont3 uni}, it follows that $u = v$ and $m=m'$.\\
%Fix $\tilde{m}$ in $B$. 
In a first case, 
assume that $\beta^2$ is not an eigenvalue of $\ca_{m}^*\ca_{m}$.
Then using the same arguments as in the proof of lemma
  \ref{proj lemma}, there is a neighborhood  $W$ %$W_{\tilde{m}}$ 
	of  $m$ %$\tilde{m}$  
	in $ B $ such that for all $s$ and $t$ in $W$,
	$\beta^2$ is not an eigenvalue of either $\ca_{s}^*\ca_{s}$ or
$\ca_{t}^*\ca_{t}$ and $\| \cp_s - \cp_t \| = O( |s-t|) $, uniformly for $s$
and $t$ in  $W$. % $W_{\tilde{m}}$.
%Assume now that $m \in W_{\tilde{m}}$.
As $\cp_{m_k'} v_k = v_k$ and $\cp_{m_k} u_k = u_k$,  we may write,
\bean \label{decomp}
\ca_{m_k} u_k -  \ca_{m_k'} v_k  = 
(\ca_{m_k}  -  \ca_{m_k'}) u_k   
- \ca_{m_k'} \cp_{m_k'} (v_k -u_k)  
 -\ca_{m_k'} ( \cp_{m_k'} - \cp_{m_k} )u_k .  
\eean
We may assume that 
$\f{{m_k}- {m_k'}}{|m_k - m_k'|}  $ 
converges to  some unit vector $\bq$ 
in $\RR^p$.
By  \ref{mder}, since $B'$ is open and $m_k$ and $m_k'$ converge to $m$,
the line segment from $m_k$ to $m_k'$ is in $B'$ for all $k$ large enough.
Let $\phi$ be in $H^1_0(R)$.
By \eqref{Am1} and \ref{mder}, if $[m_k, m_k'] \subset B'$,
\bean \label{taylor int}
\f {\ca_{m_k} \phi - \ca_{m'_k} \phi}{|m_k - m_k'|} (x)=
\int_R  \int_0^1  \nabla_m \bH_{m_k' + t(m_k - m_k')} (x,y) \cdot 
\f {m_k - m_k'}{|m_k - m_k'|}  \phi(y) d t d \sigma(y),
\eean
and since   by  \ref{mder},
\bean \label{taylor int2}
\int_0^1  \nabla_m \bH_{m_k' + t(m_k - m_k')} (x,y) \cdot 
\f {m_k - m_k'}{|m_k - m_k'|}  \ri  \nabla_m \bH_{m} (x,y) \cdot 
\bq   ,
\eean
as $k \ri \infty$, uniformly in $(x,y)$, it follows that 
$\f{\ca_{m_k} - \ca_{m'_k}}{|m_k - m_k'|}$
 converges
to $\p_\bq \ca_m $ in operator norm
and therefore
$\f{\ca_{m_k}- \ca_{m_k'}}{|m_k - m_k'|} u_k $ 
converges strongly to $\p_\bq \ca_m  u $.
%In a first case, assume that  $\beta^2$ is not an eigenvalue of  
%$\ca_{m}^* \ca_{m} $.
 $\f{\cp_{m_k'} - \cp_{m_k}}{|m_k - m_k'|} u_k$ 
is  
bounded,
%$W_{\tilde{m}}$.
%\eqref{proj_est}, 
so after possibly extracting a subsequence we may assume
by lemma \ref{intermediate lemma}
 that
$\ca_{m_k'} \f{ \cp_{m_k'} - \cp_{m_k} }{|m_k - m_k'|} u_k$ converges strongly to some
$\ca_{m} w$ for some $w$ in $H^1_0(R)$. 
Now, by \eqref{lower bound},  \eqref{contra_cont3 uni}, and \eqref{decomp}, 
$\cp_{m_k'} \f{v_k  - u_k }{|m_k - m_k'|}$ is also bounded, thus we 
 may assume that
$\ca_{m_k'} \cp_{m_k'} \f{v_k -u_k}{|m_k - m_k'|}$ converges strongly to some
$\ca_{m} z$ for some $z$ in $H^1_0(R)$.  
Altogether,  we obtain at the limit thanks to \eqref{contra_cont3 uni}
and \eqref{decomp},
\bea
\p_\bq \ca_m  u = \ca_{m} (w+z).
\eea
As $u \neq 0$, this contradicts  \ref{inter}. \\
In the second case, $\beta^2$ is  an eigenvalue of  
$\ca_{m}^* \ca_{m} $. Let $\epsilon>0$ be such that 
$\ca_{m}^* \ca_{m} $ has no eigenvalue in $(\beta^2 - 2 \epsilon, \beta^2 )$.
Now, let $\cp_t'$ be the orthogonal projection in $H^1_0(R)$ on the span
of eigenvectors of $\ca_{t}^* \ca_{t} $ corresponding to eigenvalues greater than
$\beta^2 - \epsilon$. The same argument as above may be repeated by using
 $\cp_t'$  in place of $\cp_t$. $\Box$\\
%Then for all $s$  in a neighborhood $W$ of $m$
%define $F_{s}$ to  be the subspace spanned by the eigenvectors of 
%$\ca_{s}^* \ca_s$ corresponding to eigenvalues strictly greater than $\beta^2 - \epsilon$ and
 %$R_{s}$ be the orthogonal projection in $H^1_0(R)$ on $F_{s}$.
%We can choose $W$ small enough so that $ \| R_{s}  - R_t \| = O( |s  - t|)$
%uniformly for all $s$ and $t$ in $W$.
%As $E_{s} \subset F_{s} $, the previous argument can be repeated using
 %$R_{m_k} $  in place of $\cp_{m_k}$ and $R_{m_k'} $  in place of $R_{m_k'}$.
%$\Box$\\

\begin{prop}\label{first Em prop}
Fix $\beta>0$ and define $E_m$ as above. Fix $m_0$ in $B$ and $v_0 \neq 0$ in $H^1_0(R)$.
There is a positive $C_{v_0}$ such that for all $m$ in $B$ and all $u$ in $E_m$, 
\bean \label{main est}
	 \| \ca_m u -  \ca_{m_0} v_0  \| \geq
 C_{v_0} |m - m_0|.
\eean
\end{prop}
\textbf{Proof:}
Arguing by contradiction, assume that there is a sequence $m_k$  in $B$,  and a sequence
$u_k$ in $E_{m_k}$,
	\bean \label{contra_cont3}
	 \| \ca_{m_k} u_k -  \ca_{m_0} v_0  \| < 
 \f{1}{k}  |m_0 - m_k|.
\eean
By compactness, we may assume that $m_k$ converges to some $m$ in $B$.
By \ref{U1}, $\ca_{m_0} v_0 \neq 0$. By \eqref{lower bound} and
 \eqref{contra_cont3}, the sequence 
$u_k$ is bounded, so after extracting a subsequence, we may assume that $u_k$
converges weakly to some $u$. By lemma \ref{intermediate lemma}, 
$\ca_{m_k} u_k$ converges strongly to $\ca_{m} u$. Combining \ref{U1} and
\eqref{contra_cont3}, it follows that $u = v_0$ and $m=m_0$.
Now, as $v_0$ is in $E_m$ and $m=m_0$,
 inequality \eqref{contra_cont3} contradicts \eqref{main est uni}.
$\Box$\\

\subsection{The discrete case}

\begin{lem}\label{inf norm boundedness}
There are two positive constants $C_0, C_1$ such that for all
$u$ in $H^1_0(R)$ and all $m$ in $B$,
\bea
\sup_{V} | \ca_m u | \leq C_0 \| u \|,\\
\sup_{V} |\nabla_x \ca_m u | \leq C_1 \| u \|.
\eea
\end{lem}
\textbf{Proof:}
This is clear due to \ref{H C1}.

\begin{thm} \label{main disc theorem}
%Fix two positive constants $A_1$ and $A_2$.
Fix $\beta>0$ and define $E_m $ as previously.
There is an integer $N$ such that for all $m, m'$ 
in $B$, all $u$ in $E_m$, all $v$ in $E_{m'}$,
% with
%$\|v \| \leq A_2$
%	such that if 
	%$  \sum_{j=1}^{M_{k}} C'(j,k) | 
%\ca_{m'}v (P_{j,k})|^2 \geq A_1^2$ 
and all $k >N$ in $\NN$,
\bean \label{main est disc double}
	\left(  \sum_{j=1}^{M_{k}} C'(j,k) | 
(\ca_{m}u - \ca_{m'}v) (P_{j,k})|^2 \right)^{\f12} \no \\
\geq
\f{C}{2} \left(  \sum_{j=1}^{M_{k}} C'(j,k) | 
(\ca_{m'}v) (P_{j,k})|^2 \right)^{\f12} |m - m'|,
\eean
where $C$ is the same constant as in theorem 
\ref{unif with Em}.
\end{thm}
\textbf{Proof:}
We first show that if $m_k' \in B$ and $v_k$ is a sequence such that $v_k $ is in $E_{m_k'}$,
and for a sequence $r_k \ri \infty$,
$
\ds \sum_{j=1}^{M_{r_k}} C'(j, r_k) | 
(\ca_{m_k'}v_k) (P_{j,k})|^2=1,
$
then \\
(i). $\int_V |\ca_{m_k'} v_k|^2$ converges to 1, \\
(ii). $\exists N \in \NN, \forall k >N, \|\ca_{m_k'}v_k \| \geq \f12$. \\
To prove (i).  we note that it follows from \ref{Q1} that
\bean
\label{1 minus}
| 1 - \int_V |\ca_{m_k'} v_k|^2 | = 
O(\f{1}{M_{r_k}^{\f{1}{\dim V}}})  \sup_V |\nabla_x ( \ca_{m_k'}v_k)^2 |.
\eean
Arguing by contradiction assume that a subsequence of $\int_V |\ca_{m_k'} v_k|^2$
diverges to infinity. We also denote $\int_V |\ca_{m_k'} v_k|^2$ that subsequence to 
ease notations. 
By \eqref{1 minus},
\bea
  \sup_V |\nabla_x ( \ca_{m_k'}v_k)^2 | /
\int_V |\ca_{m_k'} v_k|^2 \ri \infty,
\eea
so by \eqref{lower bound}
\bea
  \sup_V |\nabla_x ( \ca_{m_k'}v_k)^2 | /
\|  v_k \|^2 \ri \infty.
\eea
%Since for any $v$ in $H^1_0(R)$ the sup norm of
%$ \nabla_x ( \ca_{m_k'}v)^2 = 2  \ca_{m_k'}v \nabla_x \ca_{m_k'}v $
%s bounded above by a constant times $\| v \|$, this is a contradiction. \\
This contradicts  lemma \ref{inf norm boundedness}.
Thus $\int_V |\ca_{m_k'} v_k|^2$ is bounded, so by \eqref{lower bound}, lemma 
\ref{inf norm boundedness}, and \eqref{1 minus}, (i). is proved.  (ii). is then clear. \\
%\eqref{1 minus}, $\int_V |\ca_{m_k'} v_k|^2 $ converges to 1.\\
%We first show that if $v \in E_{m'}$ such that $\| v \| \leq A_2$ 
%and there is an $N$ in $\NN$ such that for all $k>N$,
%$  \sum_{j=1}^{M_{k}} C'(j,k) | 
%\ca_{m'}v (P_{j,k})|^2 \geq A_1^2$,
%then $\| \ca_{m' } v \| \geq A_1$.
%By \ref{Q1},
%\bea \sum_{j=1}^{M_{r_k}} C'(j,r_k) | 
%\ca_{m'}v (P_{j,r_k})|^2  -
%\int_V  | \ca_{m'}v|^2 
%= O(\f{1}{M_n^{\f{1}{\dim V}}})  \sup_V |\nabla_x (\ca_{m'} v)^2 |.
%\eea
%Since $\| v \| \leq A_2$, by \ref{H C1},
 %$\sup_V |\nabla_x (\ca_{m'} v)^2|$ is bounded and the claim follows by
%letting $k$ tend to infinity.\\
Since $\ca_m, \, \ca_{m'}$ are linear operators and $E_m, \, E_{m'}$ 
are linear spaces we only need to show this estimate in the case where
$ \ds \sum_{j=1}^{M_{k}} C'(j,k) | 
(\ca_{m'}v) (P_{j,k})|^2=1$. 
From \eqref{main est disc double}, arguing by contradiction, assume that 
%there is a positive $\alpha$, 
%a sequence of positive numbers $C_k$ converging to zero, 
there is a sequence  $r_k$ in $\NN$ diverging to
infinity, that
there are two sequences $m_k, m_k'$ in $B$ with $m_k \neq m_k'$
for all $k$, and a sequence
$u_k$ in $E_{m_k}$ and $v_k$ in $E_{m_k'}$
such that %$\| v_k \| \leq A_2$, $\| \ca_{m_k'} v_k \| \geq A_1$
$ \ds \sum_{j=1}^{M_{r_k}} C'(j,r_k) | 
(\ca_{m'}v_k) (P_{j,r_k})|^2=1$
and
 \bean \label{contra_cont m mprime}
	 \sum_{j=1}^{M_{r_k}} C'(j,r_k) | 
(\ca_{m_k}u_k - \ca_{m_k'}v_k ) (P_{j,r_k})|^2 <  \f{C^2}{4}|m_k - m_k'|^2 .
\eean
It follows that $ \ds \sum_{j=1}^{M_{r_k}} C'(j,r_k) | 
(\ca_{m'}u_k) (P_{j,r_k})|^2$ is bounded.
Thus by (i). and \eqref{lower bound}, $\| u_k\|$ is bounded.
Note that $\| v_k\|$ is also  bounded by  (i). and \eqref{lower bound}.
It now follows that 
$\sup_V |\nabla_x (\ca_{m_k}u_k - \ca_{m_k'}v_k)^2 |$ is 
also bounded. 
Thus
\bea
\sum_{j=1}^{M_{r_k}} C'(j,r_k) | 
(\ca_{m_k}u_k - \ca_{m_k'}v_k ) (P_{j,r_k})|^2  -
\int_R  |\ca_{m_k}u_k - \ca_{m_k'}v_k'|^2
\eea
converges to zero by assumption \ref{Q1}.
We may assume  by compactness that $m_k$ converges to some $m$ in $B$,
$m_k'$ converges to some $m'$ in $B$,
$u_k$ converges weakly to $u$ in $H^1_0(R)$, and 
$v_k$ converges weakly to $v$ in $H^1_0(R)$. 
Note that necessarily $\| \ca_{m'} v \| =1$.
As $\ca_{m_k}u_k, \ca_{m_k'}v_k$ converge strongly to $\ca_{m}u, \ca_{m'}v'$,
%by (ii). $\| \ca_{m'}v'  \| \geq \f12$.
%By 
we have found that 
 \bea
\sum_{j=1}^{M_{r_k}} C'(j,r_k) | 
(\ca_{m_k}u_k - \ca_{m_k'}v_k ) (P_{j,r_k})|^2  -
\int_R  |\ca_{m}u - \ca_{m'}v|^2
\eea
converges to zero. 
The condition $m \neq m'$ would then contradict 
%theorem 
%\eqref{main est uni} because of \eqref{contra_cont m mprime}
%and 
\ref{U1} 
since $\| \ca_{m'} v \| =1$.
After extracting a subsequence we may assume that 
$\f{ m_k - m_k'}{|m_k - m_k'|}$ converges to some $\bq $ in $\RR^p$
with $|\bq |=1$. 
Next, we want to show that 
\bean \label{to zero}
\sum_{j=1}^{M_{r_k}} C'(j,r_k) | 
(\f{\ca_{m_k}u_k - \ca_{m_k'}v_k}{|m_k - m_k'|}) (P_{j,r_k})|^2  -
\int_R  |\f{\ca_{m_k}u_k - \ca_{m_k'}v_k'}{|m_k - m_k'|}|^2
\eean 
is also convergent to zero. To that effect, we write
\bean \label{decomp2}
\f{\ca_{m_k} u_k -  \ca_{m_k'} v_k}{|m_k - m_k'|}  = 
\f{\ca_{m_k}  -  \ca_{m_k'}}{|m_k - m_k'|} u_k   
- \ca_{m_k'} \cp_{m_k'} \f{v_k -u_k}{|m_k - m_k'|}  
 -\ca_{m_k'} \f{ \cp_{m_k'} - \cp_{m_k} }{|m_k - m_k'|}u_k,  
\eean
and we first assume that $\beta^2$ is  not an eigenvalue of  
$\ca_{m}^* \ca_{m} $.
We explained in the proof of theorem \ref{unif with Em}
that $\f{\ca_{m_k} - \ca_{m'_k}}{|m_k - m_k'|}$
 converges
to $\p_\bq \ca_m $ in operator norm. 
Thanks to assumption \ref{mder} a similar argument can be carried out to show that 
$\nabla_x \f{\ca_{m_k} - \ca_{m'_k}}{|m_k - m_k'|} u_k$
 converges
to $\nabla_x \p_\bq \ca_m u $ in the sup norm over $V$,
thus $\nabla_x \f{\ca_{m_k} - \ca_{m'_k}}{|m_k - m_k'|} u_k$ and
$\f{\ca_{m_k} - \ca_{m'_k}}{|m_k - m_k'|} u_k$ 
are bounded
and by assumption \ref{Q1},
\bea
\sum_{j=1}^{M_{r_k}} C'(j,r_k) | 
(\f{\ca_{m_k}- \ca_{m_k'}}{|m_k - m_k'|}) u_k (P_{j,r_k})|^2  -
\int_R  |\f{\ca_{m_k}- \ca_{m_k'}}{|m_k - m_k'|} u_k |^2  \, \ri 0.
\eea 
Similarly, we can   argue that 
$\nabla_x \ca_{m_k'} \f{ \cp_{m_k'} - \cp_{m_k} }{|m_k - m_k'|}u_k$
and
$\ca_{m_k'} \f{ \cp_{m_k'} - \cp_{m_k} }{|m_k - m_k'|}u_k$
 are  bounded in the sup norm over $V$ since 
$\f{ \cp_{m_k'} - \cp_{m_k} }{|m_k - m_k'|}u_k$ is bounded in $H^1_0(R)$
and thus
\bea
\sum_{j=1}^{M_{r_k}} C'(j,r_k) | 
\ca_{m_k'} \f{ \cp_{m_k'} - \cp_{m_k} }{|m_k - m_k'|}u_k (P_{j,r_k})|^2  -
\int_R  |\ca_{m_k'} \f{ \cp_{m_k'} - \cp_{m_k} }{|m_k - m_k'|}u_k |^2  \, \ri 0.
\eea 
It now follows from 
\eqref{contra_cont m mprime} and 
\eqref{decomp2}
that 
\bea
\sum_{j=1}^{M_{r_k}} C'(j,r_k) | \ca_{m_k'} \cp_{m_k'} \f{v_k -u_k}{|m_k - m_k'|}  
(P_{j,r_k})|^2  
\eea
is also bounded.  We then claim by (i). and \eqref{lower bound}
that 
$\cp_{m_k'} \f{v_k -u_k}{|m_k - m_k'|}$ is bounded in $E_{m_k'}$,
so $ \ca_{m_k'} \cp_{m_k'} \f{v_k -u_k}{|m_k - m_k'|} $ and
$ \nabla_x \ca_{m_k'} \cp_{m_k'} \f{v_k -u_k}{|m_k - m_k'|}$ 
are 
bounded in sup norm.
Altogether, recalling \eqref{decomp2} and assumption \ref{Q1} , 
we have proved that  \eqref{to zero}
converges to zero.
But now, by \eqref{contra_cont m mprime}
we obtain that for $k$ large enough,
$ \|  \ca_{m_k}u_k - \ca_{m_k'} v_k \| \leq 
\f34 C | m_k - m_k'|$.
Fix $\alpha $ in $(0,1)$ and recall 
 $ \ds \sum_{j=1}^{M_{r_k}} C'(j,r_k) | 
(\ca_{m'}v_k) (P_{j,r_k})|^2=1$. For some $k$ large enough,
$1 - \alpha < \|  \ca_{m_k'} v_k \|$
thus 
$ \|  \ca_{m_k}u_k - \ca_{m_k'} v_k \| \leq 
\f34 \f{C}{1 -\alpha }  \|  \ca_{m_k'} v_k \| | m_k - m_k'|$.
Since $u_k \in E_{m_k}$, $v_k \in E_{m_k'}$,
 and $m_k \neq
m_k'$
this contradicts \eqref{main est uni} for $\alpha$ close enough to zero.\\
%we obtain at the limit
%$ \|  \ca_{m}u - \ca_{m'} v \| \leq 
%\f12 C | m - m'|$. 
%By \eqref{main est uni} this implies that  $m=m'$, and then by 
%assumption \ref{U1}, $u=v$.
In the second case, $\beta^2$ is  an eigenvalue of  
$\ca_{m}^* \ca_{m} $. As in the proof of theorem \ref{unif with Em}
we set 
$\epsilon>0$ to be such that 
$\ca_{m}^* \ca_{m} $ has no eigenvalue in $(\beta^2 - 2 \epsilon, \beta^2 )$
and we work with 
  $\cp_t'$, the orthogonal projection in $H^1_0(R)$ on the span
of eigenvectors of $\ca_{t}^* \ca_{t} $ corresponding to eigenvalues greater than
$\beta^2 - \epsilon$. The same argument as above may be repeated by using
 $\cp_t'$  in place of $\cp_t$. $\Box$\\

\begin{thm} \label{v0 disc theorem}
Fix $\beta>0$ such that \eqref{lower bound} holds.
Fix $m_0$ in $B$ and $v_0 \neq 0$ in $H^1_0(R)$.
There is  an integer $N$ such that for all $m$ in $B$, all $u$ in $E_m$, 
and  all $k >N$ in $\NN$,
\bean \label{main est disc}
	\left(\sum_{j=1}^{M_{k}} C'(j,k) | 
(\ca_{m}u - \ca_{m_0}v_0 ) (P_{j,k})|^2 \right)^{\f12}
\geq
\f{C_{v_0}}{2}  |m - m_0|,
\eean
where $C_{v_0}$ is the same constant as in proposition \ref{first Em prop}.
\end{thm}
\textbf{Proof:}
From \eqref{main est disc}, arguing by contradiction, assume that 
%there is a positive $\alpha$, 
%a sequence of positive numbers $C_k$ converging to zero, 
there is a sequence $m_k$ in $B$ with $m_k \neq m_0$ for $k \geq 1$,
 a sequence  $r_k$ in $\NN$ diverging to
infinity, and a sequence
$u_k$ in $H^1_0(R)$ such that $u_k \in E_{m_k}$ and
 \bean \label{contra_cont5}
	 \big( \sum_{j=1}^{M_{r_k}} C'(j,r_k) | 
(\ca_{m_k}u_k - \ca_{m_0}v_0 ) (P_{j,r_k})|^2 \big)^{\f12}<  \f{C_{v_0}}{2} |m_k - m_0| .
\eean
%The same argument as in the proof 
Point (i). in the proof 
of theorem \ref{main disc theorem}
 shows that 
$\| u_k\|$ is bounded.
Next, we may assume  by compactness that $m_k$ converges to some $m$ in $B$,
and 
$u_k$ converges weakly to $u$ in $H^1_0(R)$.
If $m \neq m_0$ then as $k \ri \infty$ we can contradict \eqref{main est},
%, and we can argue that $u$ is $E_m$
%since .
thus $m=m_0$ and by assumption \ref{U1} $u=v_0$.
We can then repeat the argument 
in the proof of theorem \ref{main disc theorem} to show that  this implies that for 
any $\alpha$ in $(0,1)$, for
$k$
large enough, $\| \ca_{m_k}u_k - \ca_{m_0}v_0 \| \leq \f34 \f{1}{1-\alpha} C_{v_0} |m_k - m_0| $
which contradicts \eqref{main est}
since $m_k \neq m_0$, and $u_k \in E_{m_k}$ and $v_0=u$ is in $E_m$.
%Then necessarily $u=v_0$ 
$\Box$\\

%Fix $A>0$.  Recalling the definition \eqref{setS}
%of the subspace $S$ of $H^1_0(R)$,
%fix $M_3>0$ and define (or there is)?
%\bean \label{setS prime}
%	S' = \{ \varphi \in S: \forall m \in B, \| (I- P_m) \varphi \| \leq M_3 \|P_m \varphi \|
	%\}
%\eean
%PARAGRAPH about estimating these constants based on the diameter  of B, $A_1, A_2$,
%and min of $\nabla \ca_m u $, $\| u \|=1, u \in E_m$
\subsection{Local estimate of the constant in theorem \ref{unif with Em}}
\label{constant estimates}
Let $\tilde{m}$ be in $B$. In a first case, 
assume that $\beta^2$ is not an eigenvalue of $\ca_{\tilde{m}}^*\ca_{\tilde{m}}$.
Then using the same arguments as in the proof of lemma
  \ref{proj lemma}, there is a neighborhood  $W$ %$W_{\tilde{m}}$ 
	of  $\tilde{m}$  
	in $ B $ such that for all $m$ and $m'$ in $W$,
	$\beta^2$ is not an eigenvalue of either $\ca_{m}^*\ca_{m}$ or
$\ca_{m'}^*\ca_{m'}$ and the estimate for orthogonal projections
$\| \cp_m - \cp_{m'} \| = O( |m - m'|) $ holds  uniformly for $m$
and $m'$ in  $W$. 
In fact, by the integral formula \eqref{int for proj},
 the factorization \eqref{factorization}, and \ref{mder},
it follows that $\cp_m$ has a continuous derivative in $m$ for $m$ in $W$.
By possibly shrinking $W$, we may assume that $W$ is closed and convex.
%Fix two positive constants $ A_1 $ and $ A_2$ as in theorem
 %\ref{unif with Em}.
Let 
\bea
T_{m'} = \{ v \in E_{m'} :  \| \ca_{m'} v \| = 1\}.
\eea
We then write for $u$ in $E_m$ and $v$ in $T_{m'}$, thanks 
to assumption \ref{mder},
\bean \label{write out}
&&\ca_m u - \ca_{m'} v \no \\
&=& \ca_m \cp_m (u -v) - \ca_m (\cp_{m'} - \cp_m) v - (\ca_{m'} - \ca_{m}) v \no \\
&=& \ca_m \cp_m (u -v) - \ca_m \nabla \cp_m \cdot (m' -m  )v 
- \nabla \ca_{m} \cdot ( m' -m ) v + o(|m - m'|),
\eean
where the remainder $o(|m - m'|)$ does not depend on $u$ in $E_m$ or $v$ 
in $T_{m'}$.
%Since $v \in T_{m'}$, it can also be made dependent on $A_2$ only, and otherwise
%be independent of $v$. 
\begin{prop}
\bean
\label{inf of dist}
\inf_{\bq \in \RR^p, |\bq|=1, m, m' \in W} 
{\rm dist } ((\p_\bq \ca_m  + \ca_m \p_\bq \cp_m )T_{m'}, \ca_m E_m ) >0
\eean
and 
the constant $C$
in theorem \ref{unif with Em} can be asymptotically equal to 
this inf if %$u$ and $v$ are in a
$B$ is reduced to the 
 small neighborhood $W$ 
of $\tilde{m}$.
\end{prop}
\textbf{Proof:}
We first note that the sets  
$(\p_\bq \ca_m  + \ca_m \p_\bq \cp_m )T_{m'}$ and $ \ca_m E_m$
do not intersect due to assumption \ref{inter}.
Since  $(\p_\bq \ca_m  + \ca_m \p_\bq \cp_m )T_{m'}$  is a compact set 
in $L^2(V)$ and $ \ca_m E_m$ is a finite dimensional subspace of $L^2(V)$,
it follows that for any fixed $\bq \in \RR^p$ and
fixed $m, m '$ in $W$, 
\bean \label{pos dist}
 {\rm dist } ((\p_\bq \ca_m  + \ca_m \p_\bq \cp_m )
T_{m'}, \ca_m E_m ) >0.
\eean
Arguing by contradiction, if the inf in \eqref{inf of dist} is zero,
then there is a sequence $\bq_k$ in $\RR^p$ with $|\bq_k| = 1$, two sequences
$m_k$ and $m_k'$ in $W$, a sequence $v_k$ in $T_{m'}$
and $u_k$ in $E_m$ such that
\bean \label{lim of inf}
\lim_{k \ri \infty} \|  
(\p_{\bq_k} \ca_{m_k}  + \ca_{m_k} \p_{\bq_k} \cp_{m_k} ) v_k -  \ca_{m_k} 
u_k
\| =0.
\eean
By compactness, we may assume without loss of generality that $m_k$ converges to some $m$ in $W$, $m_k'$ converges to some $m'$ in $W$, 
$\bq_k$ converges to some $\bq$ in $\RR^p$ with $|\bq |=1$, 
$v_k$ converges weakly to some $v$ in $H^1_0(R)$.
Then, we argue by
\ref{mder} and the definition of the neighborhood $W$
that
$(\p_{\bq_k} \ca_{m_k}  + \ca_{m_k} \p_{\bq_k} \cp_{m_k} ) v_k$
converges strongly to 
$(\p_{\bq} \ca_{m}  + \ca_{m} \p_{\bq} \cp_{m} ) v$.
By \eqref{lim of inf} and \eqref{lower bound}, $u_k$ is also bounded 
in $H^1_0(R)$: 
we may assume that it converges weakly to some $u$ in $H^1_0(R)$.
 Since $\cp_{m_k'} v_k = v_k$, 
as $\| \cp_{m_k'} - \cp_{m'}\| \ri 0$, and
$\cp_{m'}$ is compact, it follows that $v_k$ is strongly convergent to $v$,
so $v \in T_{m'}$. Similarly, $u \in E_m$.
As at the limit 
$(\p_{\bq} \ca_{m}  + \ca_{m} \p_{\bq} \cp_{m} ) v
 -  \ca_{m} 
u= 0$, this contradicts \eqref{pos dist}. \\
By \eqref{write out}, we see that the constant $C$
in theorem \ref{unif with Em} can be asymptotically equal to 
the inf in \eqref{inf of dist} if $u$ and $v$ are in a small neighborhood $W$ 
of $\tilde{m}$, since $\cp_m (u -v)/|m' - m|$ is in the linear space $E_m$. $\Box$\\
In the case where $\beta^2$ is  an eigenvalue of 
$\ca_{\tilde{m}}^*\ca_{\tilde{m}}$, we then use 
$\cp_m' $  in place of $\cp_m$ 
%and $\cp_{m_k'}' $  in place of $\cp_{m_k'}$,
where $\cp_m' $ %and  $\cp_{m_k'}' $ 
was defined in the proof of
proposition \ref{first Em prop} and repeat the same argument to find a local constant.\\\\
%880
\textbf{Remark:} The distance in \eqref{inf of dist} depends on $\beta$. This distance 
is increasing in $\beta$. Indeed, if $0< \beta' \leq \beta$,
then making the dependence of the space $E_m$ on $\beta$ explicit,
$E_{m, \beta} \subset E_{m, \beta'}$. Similarly 
$T_{m', \beta} \subset T_{m', \beta'}$ and thus the  distance in \eqref{inf of dist} is increasing 
in $\beta$. We will show in future work that if $\ca_m$ is defined as in 
example \ref{fracture} or \ref{sec Fault in elastic} the range of $\ca_m$
is dense, thus \eqref{inf of dist} converges to zero as $\beta$ tends to zero.

\section{Application to solving a passive  inverse elasticity problem by use 
of neural networks} \label{Application} 
\subsection{Physical and numerical interpretation of theorems \ref{main disc theorem} and \ref{v0 disc theorem}}
Fix  $\beta $ and $k$ such that formula  \eqref{main est disc double} holds. To simplify notations in this section, since $k$ is fixed,
set $M_k = M$, $P_{j,k} = P_{j}$, and $C'(j,k) = C'(j)$.
Let $S$ be the subset of $H^1_0(R) \times B$ defined by
\bea
S= \{ (u,m): u \in E_m, m \in B,   \sum_{j=1}^{M} C'(j) | 
\ca_{m}u (P_{j})|^2=1  \}.
\eea
According to
Theorem \ref{main disc theorem}
we can define a function 
\bea
\psi:  \{ (\ca_m u (P_{j}))_{1\leq j \leq M} \in \RR^{M} :
(u,m) \in S \}
\ri B\\
(\ca_m u (P_{j}))_{1\leq j \leq M} \ri m,
\eea
and $\psi$ is Lipschitz continuous.
In practice, our assumptions on the set $S$ can be interpreted as follows:
we assume that we have sufficiently many measurement points
$M$, that the magnitude of the measurements is
large enough. Normalize these measurements
using a discrete $l^2$ norm. 
%significant enough (larger in norm than $A_1$), and that there is
%a physical upper bound 
%on the forcing term ($A_2$ in the definition of $S$). 
%If these conditions are satisfied
Then $m$ can be reconstructed from
the measurements and the reconstruction is Lipschitz stable.
Another important implication of     Theorem \ref{main disc theorem}
is that 
the Lipschitz stability constant for reconstructing $m$
is inversely proportional to the magnitude of the measurements.
%the function $\psi$ can be approximated by a neural network.
Since the function $\psi$ defined above is Lipschitz regular,
it can be approximated by a neural network and 
 the growth of the depth of this neural network and of the number nodes  can be estimated given  accuracy requirements.
There are by now many papers in the neural network literature
that provide upper bounds for the size of neural networks approximating
 Lipschitz functions. For example, we  refer to  
 \cite{yarotsky2017error, shen2021neural} for estimates valid if
the ReLU (Rectified Linear Unit)  function is used for activation and 
\cite{de2021approximation} if the hyperbolic tangent function is used instead.      \\
Theorem \ref{v0 disc theorem} 
suggests what may happen if $v_0$ is not in $E_m$.
Conceivably, if $\ca_{m_0} v_0$ can be approximated by some 
$\ca_{m_0} u$ with $u$ in $E_m$, the neural network approximating $\psi$
should still be able to produce an output reasonably close to 
$m_0$
  from the input $(\ca_{m_0} v_0 (P_{j}))_{1\leq j \leq M} \in \RR^{M}$,
	and formula \eqref{main est disc} establishes the regular behavior of such an output.
	% instead, see pdf relating M1 to M@
	%In practice, 
	%if we assume again that 
%	$\|v_0 \| \leq A_2$ and$   \sum_{j=1}^{M_{k}} C'(j,k) | 
%\ca_{m_0}v_0 (P_{j,k})|^2 \geq A_1^2 $
%then the residual 
%$\ca_{m_0} (I - \cp_{m_0}) v_0$ is small, so the approximation of
%$v_0$ by $\cp_{m_0} v_0$ in $E_{m_0}$ must be good. 

%\newpage

\subsection{A numerical example}
\begin{figure}[htbp]
    \centering
      \includegraphics[scale=.4]{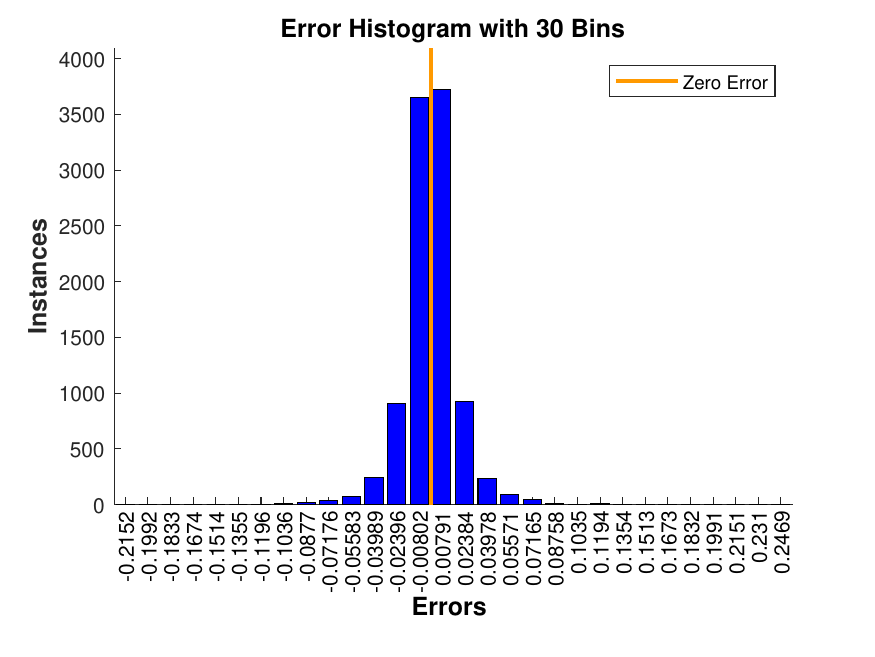}
			    \includegraphics[scale=.4]{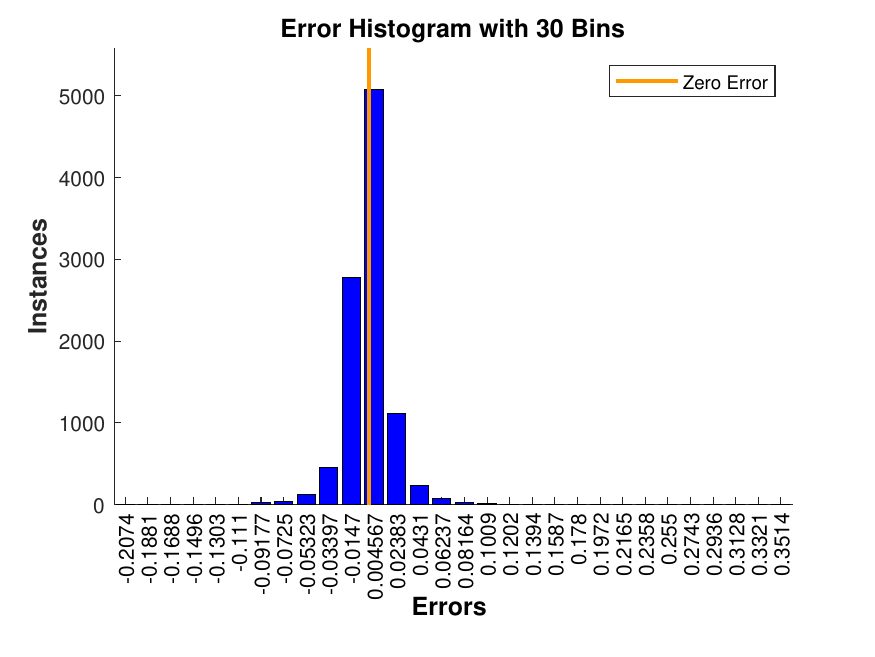}
			    \includegraphics[scale=.4]{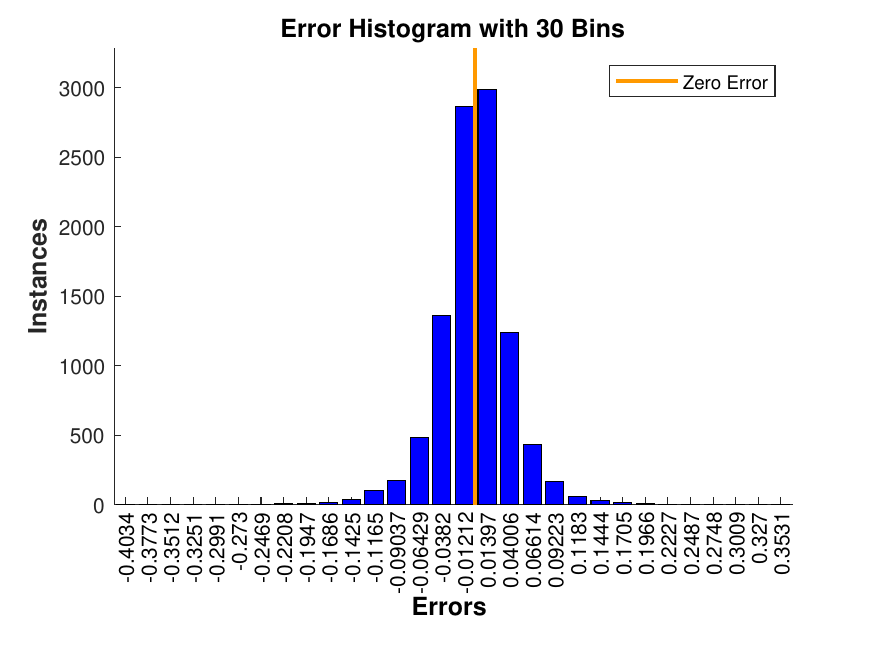}
    \caption{ Histograms of 
		errors for the learned network ${\cal N}$. 
% \\storage.wpi.edu\home\darko\My_Documents\RESEARCH\Projects for 2021\project on neural networks\121points\plotting_histograms.m
${\cal N}$ was applied to $10^4$ samples randomly selected from the
 $10^4 N_{par} $
cases used for learning. Histograms are shown for $m_1,m_2,m_3$, the three coordinates
of $m$.
$m_1,m_2,m_3$ were rescaled to normalized values
in $[0,1]$ to facilitate comparison.
}
    \label{histograms}
\end{figure}

\begin{figure}[htbp]
    \centering
      \includegraphics[scale=.4]{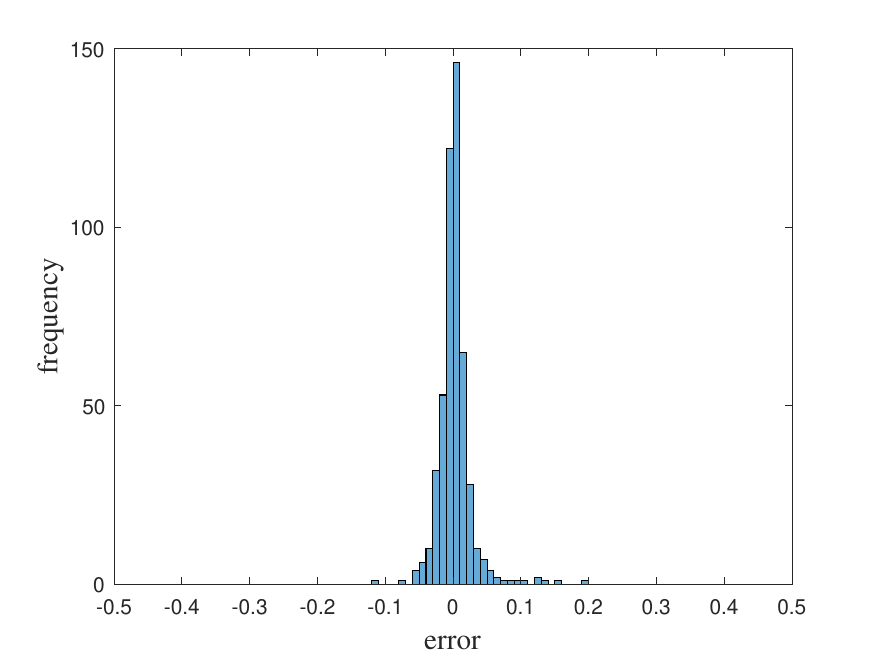}
			    \includegraphics[scale=.4]{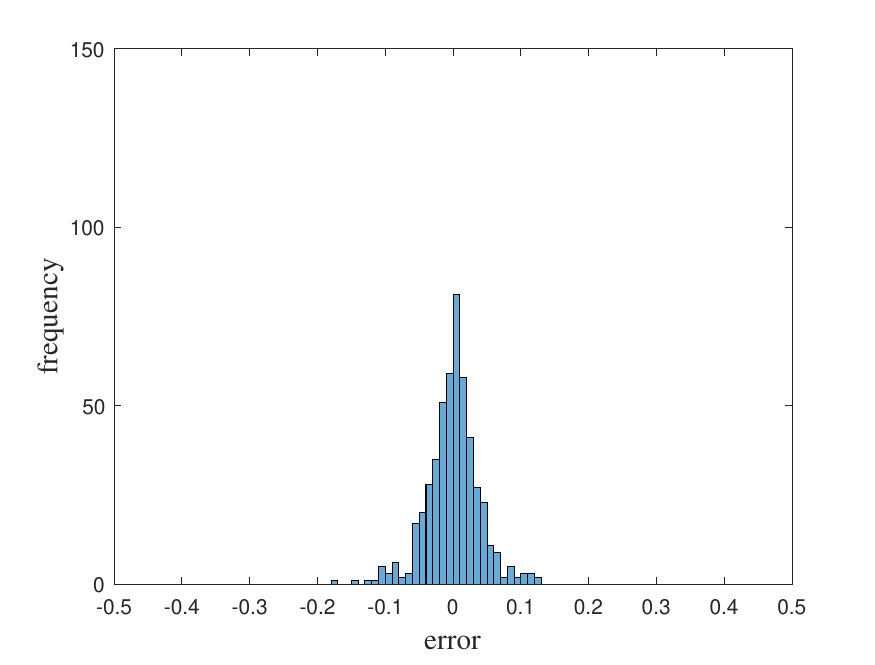}
			    \includegraphics[scale=.4]{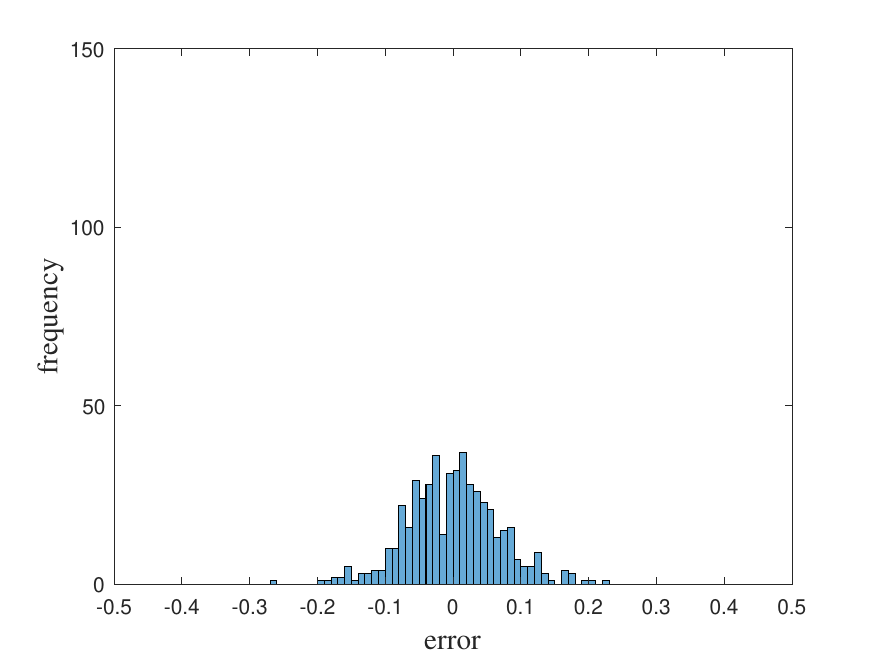}
    \caption{Errors for 
the reconstruction of $m_1$, the first coordinate of $m$, normalized to $[0,1]$ following three different methods:
by application of 
${\cal N}$, by minimization over ${\cal S}$, and minimization over ${\cal S}_0$.
These errors are for a set of $500$ draws of $m$ and random forcing vectors that were not used in the learning step. Note that the neural network 
 ${\cal N}$ reconstructs all coordinates 
	of $m$ simultaneously. Error histograms in reconstructing $m_2, m_3$ present 
	a similar profile and are not shown.}
\label{error after learning}
\end{figure}

%finding the codes: go to 
% \\storage.wpi.edu\home\darko\My_Documents\RESEARCH\Projects for 2021\project on neural networks\121points
% generate data then learn: dec21_data_for_learning.m
% verify on random selection multi_verifbestonly.m
% plotting_histograms.m plots histograms
We present a simulation 
illustrating the fault in elastic half space setting discussed in section
\ref{sec Fault in elastic} with the fault $\Gamma_m$ 
given by 
\eqref{Gammam def}, and the operator $\ca_m$ by 
\eqref{Aabd elastic}.
Here, $R$ is the square $[-150,150]^2$ in $\RR^2$ and $V$
is the square $[-200,200]$ in the plane with equation $x_3=0$.
These numbers were chosen to facilitate comparison to previous studies
\cite{volkov2020parallel, volkov2020stochastic, volkov2019stochastic, volkov2017determining}.
%\cite{volkov2017determining, volkov2019stochastic, volkov2020stochastic, volkov2020parallel}.
% sandu, IMA, invprob, JCM
% volkov2019stochastic  Darko and Sandiumenge, Joan Calafell
% volkov2017determining Pure and Applied Geophysics
% volkov2020parallel IMA Journal of Applied Mathematics
% volkov2020stochastic Journal of Computational Mathematics}
On $V$ we choose a uniform 11 by 11 grid for the points 
$P_{j}, 1\leq j \leq M$. Since the measurements are vector fields, there is a total
of 363 scalar measurements. 
% and C(j)=?
$m=(m_1,m_2,m_3)$ is confined to the box 
$B=[-2,2] \times [-2,2] \times [-10,-60]$. Here too, these numbers relate to a
wide range of possibilities in geophysical applications 
\cite{volkov2017determining} where the length scale for $x$ is a kilometer. 
In order to achieve maximum expediency of our numerical codes, instead of fixing a threshold $\beta$ we choose $E_m$ to be the space spanned by the 
first $q=5$ singular vectors $u_i$, $i=1,..,q$, of $\ca_m$.\\
% figure with singular vectors?
The computations were done on a parallel platform using $N_{par}=20$ 
processors.
We first generated data by sampling $10^4 N_{par} $ random points $m$
in $B$. For each of these random points, we generated a realization $w$ of a vector Gaussian 
in $\RR^q$ with zero mean and identity covariance
and we formed the vector
 $(\sum_{i=1}^q w_i \ca_m u_i(P_{j}) \cdot e_l), 1\leq j \leq M, 1 \leq  l
 \leq 3$,
in $\RR^{3 M}$, where $e_l, l=1,2,3$ is the natural basis
of $\RR^3$. Finally, this vector was normalized and used as input for 
learning $m$: denote ${\cal S}$ the resulting set of $10^4 N_{par}$ 
samples in $\RR^{3 M}$.   
 Although $m$ is in $\RR^3$ and $u$ is in a $q$ -dimensional space
with $q=5$, the unknown for this problem is not embedded in $\RR^8$,
it is an 8-dimensional manifold embedded in $H^1_0(R) \times B$.
A single layer of neural networks proved to be inadequate due to the complexity
of the problem.
The learning was done on a network with three hidden layers with dimension
%size_net=[250,100,30];
$250 \times 100 \times 30$. 
%  To use a hyperbolic tangent activation for deep learning, use the tanhLayer function or the dlarray method tanh.
%transferFcn: 'tansig'
We used the hyperbolic tangent function for the activation function.
This architecture requires determining  119223 weights,
so a classic Levenberg-Marquardt backpropagation algorithm is inadequate. 
After trying several training functions available in Matlab, we found out that the
scaled conjugate gradient backpropagation algorithm
\cite{moller1993scaled} led to the best results 
while still completing the learning process in just a few hours on our parallel platform. 
Some regularization of the weights was necessary to improve the 
generalization capabilities of the network.
Let $\gamma$ be in $[0,1]$. Define a convex combination between
the mean of the square of the weights and the mean square error (MSE)
for the network with coefficients $\gamma$ and $1 -\gamma$. 
We found that setting $\gamma$ close to .2 led to best 
generalization performance for this network.\\
After 6000 back propagation steps, there was close to no measurable gain in MSE.
Let ${\cal N} $ be the resulting neural network.
We show in Figure \ref{histograms} histograms of errors for 
% \\storage.wpi.edu\home\darko\My_Documents\RESEARCH\Projects for 2021\project on neural networks\121points\plotting_histograms.m
 $10^4$ samples randomly selected from the $10^4 N_{par} $
cases used for learning. $m_1,m_2,m_3$ were rescaled to normalized values
in $[0,1]$ to facilitate comparison.
Next, we evaluated how the network performs on new data.
500 random points $m$ were drawn in $B$ and for each $m$
 a random vector 
$(\sum_{i=1}^q w_i \ca_m u_i(P_{j}) \cdot e_l), 1\leq j \leq M, 1 \leq l \leq 3$,
in $\RR^{3 M}$ was formed then normalized
to obtain a test set ${\cal T}_q$.
In Figure  \ref{error after learning} we show errors for 
the reconstruction of $m_1$ normalized to $[0,1]$ following three different methods:
\begin{enumerate}
  \item by applying the neural network  ${\cal N}$,
	 \item by minimizing over 
	  the  set ${\cal S}$ of $10^4 N_{par}$ samples, 
	 \item by minimizing over a much smaller subset ${\cal S}_0$ 
	of ${\cal S}$
	 with $10^2 N_{par}$ elements randomly selected from ${\cal S}$.
	\end{enumerate}
	Note that the neural network  ${\cal N}$ reconstructs all coordinates 
	of $m$ simultaneously. Error histograms for reconstructing $m_2, m_3$ present 
	a similar profile and are not shown.
 In Figure  \ref{accuracy an run time} we show a table comparing accuracy and run time between these three methods.
Interestingly, we observe that applying the neural network  ${\cal N}$ 
is about 1000 times faster than minimizing over ${\cal S}$ even though
applying ${\cal N}$ is about twice as accurate. 
If we use the reduced sample set ${\cal S}_0$, the minimization step is drastically
faster, but not as fast as  applying ${\cal N}$ and not nearly as accurate. \\
Next, we show how ${\cal N}$ performs 
if the input comes from some $\ca_{m_0} v_0$ plus noise where 
$v_0 $ is in $H^1_0(R) $ but not necessarily in $E_{m_0}$.
Fix $q'=50$.
 We generated  $500$ random points $m$
in $B$.
For each of these random points, we generated a realization $w$ of a vector Gaussian 
in $\RR^{q'}$ with zero mean and identity covariance
and we formed the vector
 $(\sum_{i=1}^{q'} w_i \ca_m u_i(P_{j}) \cdot e_l), 1\leq j \leq M, 1 \leq  l
 \leq 3$
in $\RR^{3 M}$.
We then added  noise to this vector by first computing its sup norm  $sn$
and adding a random vector sampled from a Gaussian distribution 
in  $\RR^{3 M}$ with zero mean and covariance given by the identity times
$sn/20$. 
% sketch?
The data was then normalized to obtain a test set ${\cal T}_{q',noisy}$.
Let ${\cal T}_{q',0}$ be the corresponding noise free test set.
 In table  \ref{accuracy2} we compare accuracy 
 for the  three methods described above applied to the test sets 
${\cal T}_{q',0}$ and  ${\cal T}_{q',noisy}$. For the test set ${\cal T}_{q',0}$
we observe that there is no significant loss of accuracy compared to the
 accuracy for the first test set  ${\cal T}$. 
Accuracy deteriorates for  ${\cal T}_{q',noisy}$, but only when applying the
neural network ${\cal N}$. In this case, applying ${\cal N}$ or minimizing over
${\cal S}$ leads to nearly identical accuracy. However, running ${\cal N}$ 
is still about 1000 times faster. 

%  cal S and cal S ' same as previously
 
% used net.performParam.normalization = 'percent';% or 'standard';

\begin{figure}
\begin{center}
\begin{tabular}{ |c|c|c|c| } 
 \hline
             & average error & load time &run time \\ 
  \hline
	${\cal N}$  & 0.015  & 0.018  & 0.049 \\ 
  \hline
	 ${\cal S}$ &0.028 & 3.5  & 47\\ 
 \hline
 ${\cal S}_0$& 0.055 & 0.053 & 0.19 \\ 
 \hline
\end{tabular}
\end{center}
  \caption{For the same  500 instances 
	of the data for the inverse problem as in figure
	\ref{error after learning}: 
	column 1, average error (absolute value) for  $m_1$ normalized to $[0,1]$
		following three different methods:
by application of ${\cal N}$, by minimization over ${\cal S}$, and minimization over ${\cal S}_0$. Column 2: load time of file containing  ${\cal N}$, 
 ${\cal S}$,  
${\cal S}_0$, in seconds. Column 3: cumulative run time for these 500 instances, in seconds.}
\label{accuracy an run time}
% errors were rescaled m_1 only
% 500 draws
\end{figure}

\begin{figure}
\begin{center}
\begin{tabular}{ |c|c|c| } 
 \hline
             & average error for ${\cal T}_{q',0}$ &  
						average error for ${\cal T}_{q',noisy}$  \\
  \hline
	${\cal N}$  & 0.0178 & 0.0286\\ 
  \hline
	 ${\cal S}$ &0.0277  & 0.0280 \\ 
 \hline
 ${\cal S}_0$& 0.0491& 0.0488\\ 
 \hline
\end{tabular}
\end{center}
% data for this table
%${\cal T}$ running_time50_nonoise
% ${\cal T}'$ running_time50
% code plotting_histograms
% data files generated by 
% multi_verifbestonly
  \caption{Column 1: average error over the test set ${\cal T}_{q',0}$
	%and ${\cal T}_{q',noisy}$
	for 
the reconstruction of $m_1$ normalized to $[0,1]$ following three  methods:
by application of ${\cal N}$,  by minimization over ${\cal S}$, and minimization over ${\cal S}_0$. 
Column 2: same quantities for the test set   ${\cal T}_{q',noisy}$.}
\label{accuracy2}
% errors were rescaled m_1 only
% 500 draws
\end{figure}

\newpage
\vskip 10pt
\Large{\bf{Funding}} \\
\normalsize
This work was supported by
 Simons Foundation Collaboration Grant [351025].
%\vskip 10pt

\bibliography{ref}{}
\bibliographystyle{abbrv}
%\bibliographystyle{apalike}
 
%\bibitem{S} E. P. Stephan, A boundary integral equation method for three-dimensional crack problems in elasticity, Mathematical Methods In The Applied Sciences, 1986 ,Volume: 8 Issue: 4.

\end{document}